\UseRawInputEncoding
\documentclass[12pt]{amsart}
\usepackage[top=30truemm,bottom=30truemm,left=25truemm,right=25truemm]{geometry}
\usepackage{mathrsfs}

\usepackage{color}
\usepackage{bm}
\usepackage{amsfonts,amssymb}
\usepackage{dsfont}
\usepackage{amscd}
\usepackage{extarrows}
\usepackage{amsmath}
\usepackage{mathrsfs}
\usepackage{enumerate}
\usepackage{amscd}
\usepackage[all]{xy}
\usepackage{hyperref}
\numberwithin{equation}{section}
\theoremstyle{plain} 

\theoremstyle{definition} 

\newtheorem{thm}{Theorem}[section]
\newtheorem{cor}[thm]{Corollary}
\newtheorem{lem}[thm]{Lemma}

\theoremstyle{definition}
\newtheorem{defn}{Definition}[section]

\theoremstyle{remark}

\newcommand{\be}{\begin{equation}}
	\newcommand{\ee}{\end{equation}}
\newcommand{\bea}{\begin{eqnarray}}
	\newcommand{\eea}{\end{eqnarray}}
\newcommand{\ben}{\begin{eqnarray*}}
	\newcommand{\een}{\end{eqnarray*}}
\newcommand{\bt}{\begin{split}}
	\newcommand{\et}{\end{split}}
\newcommand{\bet}{\begin{equation}}
	\newcommand{\mc}{\mathbb{C}}
	\newcommand{\mr}{\mathbb{R}}
\newcommand{\mb}{\mathbb{B}^n}
	\newcommand{\ra}{\rightarrow}
	\newcommand{\beq}{\begin{equation*}}
		\newcommand{\eeq}{\end{equation*}}
	\newcommand{\bi}{\begin{itemize}}
		\newcommand{\ei}{\end{itemize}}

\newcommand{\rw}{\rightarrow}
\newcommand{\rwo}{\mapsto}

\DeclareMathOperator{\p}{P}
\DeclareMathOperator{\s}{S}

\begin{document}

\title[Strict curvature positivity of direct image bundles]
{Curvature strict positivity of direct image bundles associated to pseudoconvex families of domains}
		
\author[F. Deng]{Fusheng Deng}
\address{Fusheng Deng: \ School of Mathematical Sciences, University of Chinese Academy of Sciences\\ Beijing 100049, P. R. China}
\email{fshdeng@ucas.ac.cn}

\author[J. Hu]{Jinjin Hu}
\address{Jinjin Hu: \ School of Mathematical Sciences, University of Chinese Academy of Sciences\\ Beijing 100049, P. R. China}
\email{hujinjin21@mails.ucas.ac.cn}

\author[X. Qin]{Xiangsen Qin}
\address{Xiangsen Qin: \ School of Mathematical Sciences, University of Chinese Academy of Sciences\\ Beijing 100049, P. R. China}
\email{qinxiangsen@amss.ac.cn}

\begin{abstract}
We consider the curvature strict positivity of the direct image bundle associated to a pseudoconvex family of bounded domains.
The main result is that the curvature of the direct image bundle associated to a strictly pseudoconvex family
of bounded circular domains or Reinhardut domains are strictly positive in the sense of Nakano, even if the weight functions are not strictly plurisubharmonic. This result gives a new geometric insight about the property of strict pseudoconvexity, and has some applications in complex analysis
and convex analysis. We investigate that the main result implies a remarkable result of Berndtsson which states that, for an ample vector bundle
$E$ over a compact complex manifold $X$ and any $k\geq 0$, the bundle $S^kE\otimes\det E$ admits a Hermitian metric whose curvature is strictly positive
in the sense of Nakano, where $S^kE$ is the $k$-th symmetric product of $E$.
The two main ingredients in the argument of the main theorems are Berndtsson's estimate of the lower bound of curvature of direct image bundles and Deng-Ning-Wang-Zhou's characterization of the curvature Nakano positivity of Hermitian vector bundles in terms of $L^2$-estimate of $\bar\partial$.
 \end{abstract}

		\maketitle
		
\tableofcontents
\section{Introduction}
Let $U\subset \mc^n$ and $D\subset \mc^m$ be pseudoconvex bounded domains,
and let $\varphi$ be a smooth plurisubharmonic function
defined on some (open) neighborhood of the closure of $\Omega:=U\times D$ in $\mc^n\times\mc^m$.
For $t\in U$, we define the Hilbert space
$$E_t=\{f\in \mathcal O(D); \|f\|^2_t:=\int_{D}|f|^2e^{-\varphi_t}d\lambda_z<\infty\},$$
where $\mathcal O(D)$ is the space of holomorphic functions on $D$, $\varphi_t(z)=\varphi(t,z)$,  and $d\lambda_z$ is the Lebesgue measure on $\mc^m$.
When $t$ varies in $U$, $E_t$ is invariant as a vector space, but the inner product defined by the above norm varies if $\varphi$ is not constant with respect to $t$.
Set $E=\cup_{t\in U}E_t$ and take $\pi:E\ra U$ by setting $\pi(E_t)=\{t\}$,
then $E$ is a holomorphic vector bundle (of infinite rank) over $U$ with a Hermitian metric $h$ given by
$$h_t(f,g)=\int_D f\bar g e^{-\varphi_t}d\lambda_z,\ f,g\in E_t.$$
A fundamental result of Berndtsson is as follows.

\begin{thm}[{\cite[Theorem 1.1]{Ber09}}]\label{thm:Bern direct image bd domain}
With the above notations and assumptions, the curvature of the Hermitian vector bundle $(E,h)$ is semi-positive in the sense of Nakano,
and is strictly positive in the sense of Nakano if $\varphi$ is strictly plurisubharmnic.
\end{thm}

Our main purpose is to study strict positivity of curvature of direct image bundles defined in a similar way.
In Theorem \ref{thm:Bern direct image bd domain}, we see that the strict positivity of the curvature of $(E,h)$ comes from the strict
plurisubharmonicity of the weight function $\varphi$, which can be viewed as the strict curvature positivity of the trivial line bundle over $\Omega$
with Hermitian metric given by $e^{-\varphi}$.
In the present work, we show that the strict positivity of the curvature can come from a completely different source, namely,
the strict pseudoconvexity of the total space of the family of domains.

To state the main result, we first introduce some notions and notations.
Denote by $p:\mc^n\times\mc^m\ra\mc^n$ the natural projection, and for a set $A\subset \mc^n\times\mc^m$,
we denote $p^{-1}(t)\cap A$ by $A_t$, which is called the fiber of $A$ over $t$.
Of course we can view $A_t$ as a family of subsets in $\mc^m$ depending on the parameter $t$.


\begin{defn}\label{def:strict p.s.c family}\
\bi
\item[(1)] A \emph{family of domains} of dimension $m$ over a domain $U\subset\mc^n$ is a domain $\Omega\subset U\times\mc^m$ such that $p(\Omega)=U$;
such a family is called a family of bounded domains if all fibers $\Omega_t\subset\mc^m\ (t\in U)$ are bounded.
\item[(2)] A family of domains $\Omega$ over $U$ has $C^k$ ($k\geq 1$) boundary if  there exists a $C^k$ function $\rho(t,z)$
defined on some neighborhood of the closure $\overline{\Omega}$ of $\Omega$ in $U\times\mc^m$, such that $\Omega=\{(t,z)\in U\times\mc^n|\rho(t,z)<0\}$ and $d(\rho|_{\partial\Omega_t})\neq 0$ for all $t\in U$.
Such a function $\rho$ is  called a \emph{defining function} of $\Omega$.
\item[(3)] A family of domains $\Omega\subset U\times\mc^m$ over $U$ is called\emph{ pseudoconvx } if $\Omega$ is a pseudoconvex domain.
\item[(4)]A family of domains $\Omega\subset U\times\mc^m$ with $C^2$-boundary is said to \emph{have plurisubharmonic defining function} if it admits a defining function that is plurisubharmonic on some neighborhood of $\overline{\Omega}$ in $U\times\mc^m$, and is called \emph{strictly pseudoconvex} if it admits a defining function that is strictly plurisubharmonic on some neighborhood of $\overline{\Omega}$ in $U\times\mc^m$.
\ei
\end{defn}

If $U$ is pseudoconvex and $\Omega$ has a plurisubharmonic defining function, then $\Omega$ is a pseudoconvex family of domains.
But the opposite is not true in general.

We consider pseudoconvex families of bounded domains with certain symmetries, namely,
pseudoconvex families of bounded domains whose fibers are circular domains or Reinhardt domains.

Recall that a domain $D\subset \mc^m$ is called a circular domain if it is invariant under the action of $S^1$
on $\mc^m$ given by
$$e^{i\theta}\cdot(z_1,\cdots, z_m)=(e^{i\theta}z_1,\cdots, e^{i\theta}z_m),\ \theta\in\mr,$$
and is called a Reinhardt domain if it is invariant under the action of the torus group $T^m$
on $\mc^m$ given by
$$(e^{i\theta_1},\cdots, e^{i\theta_m})\cdot(z_1,\cdots, z_m)=(e^{i\theta_1}z_1,\cdots, e^{i\theta_m}z_m),\ \theta_i\in\mr.$$

\begin{thm}\label{thm:curvartue positive cicular domains}
Let $\Omega\subset U\times\mc^m$ be a strictly pseudoconvex family of bounded domains over $U\subset\mc^n$ and $\varphi$ be a $C^2$ plurisubharmonic function defined on some neighborhood of $\overline\Omega$ in $U\times\mc^m$.
We assume that all fibers $\Omega_t\ (t\in U)$ are (connected) circular domains in $\mc^m$ containing the origin and $\varphi(t,z)$ is $S^1$- invariant with respect to $z$.
Let $k\geq 0$ and $E^k_t$ be the space of homogenous polynomials on $\mc^m$ of degree $k$,  with inner product $h_t$ given by
$$h_t(f,g)=\int_{\Omega_t}f\bar ge^{-\varphi_t}d\lambda_z,\ f, g\in E_t^k.$$
We set $E^k=\cup_{t\in U}E^k_t$ and view it as a (trivial) holomorphic vector bundle over $U$ in the natural way.
Then the curvature of the holomorphic Hermitian vector bundle $(E^k,h)$ is strictly positive in the sense of Nakano.
\end{thm}

Similar results holds for a strictly pseudoconvex family of Reinhardt domains.

\begin{thm}\label{thm:curvartue positive Reinhardt domains}
Let $\Omega\subset U\times\mc^m$ be a strictly pseudoconvex family of bounded domains over $U\subset\mc^n$ and $\varphi$ be a $C^2$ plurisubharmonic function defined on some neighborhood of $\overline\Omega$ in $U\times\mc^m$.
We assume that all fibers $\Omega_t\ (t\in U)$ are (connected) Reinhardt domains in $\mc^m$ and $\varphi(t,z)$ is $T^m$- invariant with respect to $z$.
Then for any nonnegative integers $k_1,\cdots, k_m$, the function $\psi(t)$ defined by
$$e^{-\psi(t)}=\int_{\Omega_t}|z_1^{k_1}\cdots z_m^{k_m}|^2e^{-\varphi_t}d\lambda_z$$
is a strictly plurisubharmonic function on $U$.
\end{thm}

In fact, in Theorem \ref{thm:curvartue positive Reinhardt domains}, if we assume all fibers $\Omega_t$  has no intersection with any coordinate axis,
then $k_1,\cdots, k_m$ can be taken to be any integers (not necessarily nonnegative).

We now discuss the relation of Theorem \ref{thm:curvartue positive cicular domains} and Theorem \ref{thm:curvartue positive Reinhardt domains} with Theorem \ref{thm:Bern direct image bd domain}.
If we assume that both $\Omega$ is a product domain as in Theorem \ref{thm:Bern direct image bd domain} and $\varphi$ is strictly plurisubharmonic,
then, as observed in \cite{DZZ17}, the conclusions in Theorem \ref{thm:curvartue positive cicular domains} and Theorem \ref{thm:curvartue positive Reinhardt domains} can be deduced
from Theorem \ref{thm:Bern direct image bd domain}, with the help of some basic group representation theory.
On the other hand, as we will see in the proofs, if one of the above two assumptions is dropped,
then Theorem \ref{thm:curvartue positive cicular domains} and Theorem \ref{thm:curvartue positive Reinhardt domains} essentially go beyond Theorem \ref{thm:Bern direct image bd domain}.
Here the key point we want to emphasize about Theorem \ref{thm:curvartue positive cicular domains} and Theorem \ref{thm:curvartue positive Reinhardt domains}
is that strict pseudoconvexity of the family $\Omega$ encodes the curvature strict positivity of the (character) direct image bundles.

It is also interesting to compare Theorem \ref{thm:curvartue positive cicular domains} and Theorem \ref{thm:curvartue positive Reinhardt domains} with the main result in \cite{Ber11},
where Berndtsson shows that the curvature of the direct image bundle of the relative canonical bundle twisted with a Hermitian line bundle
associated to a K\"ahler family of compact manifolds is strictly positive in the sense of Nakano,
provided that the related Kodaira-Spencer map is nondegenerate and the curvature of the line bundle is strictly positive along fibers.
Berndtsson also gives counterexamples to this result if one of the two conditions is removed.
In connection to Berndtsson's result, one may imagine from Theorem \ref{thm:curvartue positive cicular domains} and Theorem \ref{thm:curvartue positive Reinhardt domains}
that strict pseudoconvexity of the family $\Omega$ implicitly implies nontrivial deformation of the fibers and certain curvature positivity along fibers.
From this point of view, it seems that Theorem \ref{thm:curvartue positive cicular domains} and Theorem \ref{thm:curvartue positive Reinhardt domains} provide a very deep new geometric insight about strict psedoconvexity in complex analysis and complex geometry.
It seems that more profound potential relations of Theorem \ref{thm:curvartue positive cicular domains} and Theorem \ref{thm:curvartue positive Reinhardt domains}
with Berndtsson's result deserves further study.

On the other hand, we conjecture that certain appropriate form of the converse of Theorem \ref{thm:curvartue positive cicular domains} ( or Theorem \ref{thm:curvartue positive Reinhardt domains}) holds, namely, the curvature strict positivity of the direct images implies the strict pseudoconvexity of the family $\Omega$.

We want to point out that the symmetry involved in Theorem \ref{thm:curvartue positive cicular domains} and Theorem \ref{thm:curvartue positive Reinhardt domains}
does not play essential role, and it is mainly used to avoid considering the whole space of $L^2$-holomorphic functions on $\Omega_t$ as in Theorem \ref{thm:Bern direct image bd domain} and bundles of infinite rank without local trivialization.
As mentioned above,  the key role is played by the strict psedoconvexity of the family $\Omega$.
On the other hand, as we will see, the $S^1$-symmetry
is indispensable when we apply  Theorem \ref{thm:curvartue positive cicular domains} to the study of ample vector bundles over projective manifolds.

We now discuss some consequences of Theorem \ref{thm:curvartue positive cicular domains} and Theorem \ref{thm:curvartue positive Reinhardt domains}.

\begin{cor}\label{cor:Bergman kernel strict psh}
Let $\Omega\subset U\times\mc^m$ be a strictly pseudoconvex family of bounded domains over $U\subset\mc^n$ and $\varphi$ be a $C^2$ plurisubharmonic function defined on some neighborhood of $\overline\Omega$ in $U\times\mc^m$
that satisfy the conditions in  Theorem \ref{thm:curvartue positive cicular domains} or Theorem \ref{thm:curvartue positive Reinhardt domains}.
For $t\in U$, let $K(t,z)$ be the weighted Bergman kernel of $\Omega_t$ with weight $\varphi_t$.
Then $\ln K(t,z)$ is a strictly plurisubharmonic function on $\Omega$.
\end{cor}

We believe that Corollary \ref{cor:Bergman kernel strict psh} holds for an arbitrary strictly pseudoconvex family of bounded domains, without symmetry.
But we will not discuss this topic further in the present work.

\begin{cor}\label{cor:strict minimum principle}
Let $\Omega\subset U\times\mc^m$ be a strictly pseudoconvex family of domains over $U\subset\mc^n$ and $\varphi$ be a $C^2$ plurisubharmonic function defined on some neighborhood of $\overline\Omega$ in $U\times\mc^m$.
\bi
\item[(1)] If $\Omega$ and $\varphi$ satisfy the conditions in  Theorem \ref{thm:curvartue positive cicular domains} or Theorem \ref{thm:curvartue positive Reinhardt domains}, then the function $\tilde\varphi$ defined by
$$e^{-\tilde\varphi(t)}=\int_{\Omega_t}e^{-\varphi(t,z)}d\lambda_z$$
is a strictly plurisubharmonic function on $U$.
\item[(2)]If all fibers $\Omega_t$ are tube domains of the form $X_t+i\mr^m$ with $X_t$ bounded, and $\varphi(t,z)$ does not depend on the imaginary part of $z$,
then the function $\tilde\varphi$ defined by
$$e^{-\tilde\varphi(t)}=\int_{X_t}e^{-\varphi(t,\text{Re}z)}d\lambda_{\text{Re}z}$$
is a strictly plurisubharmonic function on $U$.
\ei
\end{cor}

Taking $\varphi\equiv 0$ in Corollary \ref{cor:strict minimum principle}, we get

\begin{cor}\label{cor:volue vairation strict psh}
Let $\Omega\subset U\times\mc^m$ be a strictly pseudoconvex family of domains over $U$.
\bi
\item[(1)] If $\Omega$ satisfies the conditions in  Theorem \ref{thm:curvartue positive cicular domains} or Theorem \ref{thm:curvartue positive Reinhardt domains}, then the function given by
$$t\mapsto -\ln|\Omega_t|$$
is a strictly plurisubharmonic function on $U$,
where $|\Omega_t|$ is the Lebesgue measure of $\Omega_t\subset\mc^m$.
\item[(2)] If all fibers $\Omega_t$ are tube domains of the form $X_t+i\mr^m$ with $X_t$ bounded,
then the function given by
$$t\mapsto -\ln|X_t|$$
is a strictly plurisubharmonic function on $U$,
where $|X_t|$ is the Lebesgue measure of $X_t\subset\mr^m$.
\ei
\end{cor}

The plurisubharmonicity of $\ln K(t,z)$ in Corollary \ref{cor:Bergman kernel strict psh} was proved in \cite{Ber06},
and the plurisubharmonicity of the functions considered in Corollary \ref{cor:strict minimum principle} and Corollary \ref{cor:volue vairation strict psh}
were proved in \cite{Ber98}.
The contribution here is on the strict plurisubharmonicity of those functions.

We will explain that the above corollaries imply some parallel results in convex analysis,
following a general principle given in \cite{DJQ22} that connecting convex analysis and complex analysis.

In a similar way as in Definition \ref{def:strict p.s.c family}, we define a strictly convex family of domains in $\mr^m$ as follows.
Let $U_0\subset\mr^n$ be a domain.
By definition, a \emph{strictly convex family of domains} over $U_0$ is a convex domain $D\subset U_0\times\mr^m$ such that $p_0(D)=U_0$,
and there exists a $C^2$ strictly convex function $\rho_0(t, x)$ on some neighborhood $\tilde D$ of $\overline D$ in $U_0\times\mr^m$ such that
$$D=\{(t,x)\in \tilde D|\rho_0(t, x)<0\}$$
and $d(\rho_0|_{\partial D_{t}})\neq 0$
 for all $t\in U_0$,
where $p_0:\mr^n\times\mr^m\ra\mr^n$ is the natural projection and $D_{t}=p_0^{-1}(t)\cap D$.
The function $\rho_0$ is called a defining function of $D$.
Here a $C^2$ function is called strictly convex if its Hessian is positively definite everywhere.


\begin{cor}\label{cor:strict prekopa}
Let $D\subset U_0\times\mr^m$ be a strictly convex family of bounded domains over a domain $U_0\subset\mr^n$
and $\varphi$ be a $C^2$ convex function defined on some neighborhood of the closure of $D$ in $U_0\times\mr^m$.
Then the function $\tilde\varphi$ defined by
$$e^{-\tilde\varphi(t)}=\int_{D_t}e^{-\varphi(t,x)}d\lambda_x$$
is a strictly convex function on $U_0$.
\end{cor}

The convexity of $\tilde\varphi$ in Corollary \ref{cor:strict prekopa} can be deduced from the Pr\'ekopa's theorem \cite{Pre73},
but here we are interested in the strict convexity of $\tilde\varphi$.

Taking $\varphi\equiv 0$ in Corollary \ref{cor:strict prekopa}, we get a stronger form of the classical Brunn-Minkowski inequality in convex analysis.

\begin{cor}\label{cor:volue vairation strict convex}
Let $D\subset U_0\times\mr^m$ be a strictly convex family of bounded domains over $U_0\subset\mr^n$, then the function given by
$$t\mapsto -\ln|D_t|$$
is a strictly convex function on $U_0$.
\end{cor}

Theorem \ref{thm:curvartue positive cicular domains} has a direct application to vector bundles.
Let $\pi:E\ra X$ be a holomorphic vector bundle of rank $m$ over a complex manifold $X$.
By definition, a smooth Finsler metric on $E$ is a continuous function $h:E\ra \mr$ such that $h\geq 0$,
$h(\lambda v)=|\lambda|h(v)$ for $\lambda\in\mc$ and $v\in E$, and $h$ is smooth on $E\backslash Z_E$, where $Z_E\subset E$ is the zero section of $E$.
We call $(E,h)$ is \emph{strictly negatively curved} if $\ln h$ is strictly plurisubharmonic on $E\backslash Z_E$.
(Note that if $h$ is a smooth Hermitian metric, then $(E,h)$ is strictly negatively curved if and only if its curvature is
strictly negative in the sense of Griffiths.)

Given a smooth Finsler metric $h$ on $E$, we can define an induced Hermitian metric $\text{det}h$ on the determinant line bundle $\text{det}E=\Lambda^mE$ of $E$
via the measure $\mu$ on $E_t$ with $\mu(\{v\in E_t|h(v)\leq 1\})=1$ (see \S \ref{sec: corollary} for details), $t\in X$.
(In fact, the definition still works even if $h$ is just a singular Finsler metric.)

%

\begin{cor}\label{cor: detE strict negative}
Let $\pi:E\ra X$ be a holomorphic vector bundle over a complex manifold $X$ equipped with a smooth Finsler metric $h$.
If $(E,h)$ is strictly negatively curved, then the curvature of the induced Hermitian metric $\text{det}h$ on $\text{det}E$ is strictly negative.
\end{cor}

Motivated by the study of ample vector bundles,
we also establish a result about the curvature strict positivity of invariant direct images from another perspective.

\begin{thm}\label{thm:curvartue positive cicular domains from weight}
Let $\Omega\subset U\times\mc^m$ be a family of bounded domains over $U$ that admits a plurisubharmonic defining function,
and $\varphi$ be a $C^2$ plurisubharmonic function defined on some neighborhood of $\overline\Omega$ in $U\times\mc^m$.
We assume that all fibers $\Omega_t\ (t\in U)$ are (connected) circular domains in $\mc^m$ containing the origin and $\varphi(t,z)$ is $S^1$ invariant with respect to $z$.
Let $k\geq 0$ and $E^k_t$ be the space of homogenous polynomials on $\mc^m$ of degree $k$, with inner product $h_t$ given by
$$h_t(f,g)=\int_{\Omega_t}f\bar ge^{-\varphi_t}d\lambda_z,\ f, g\in E_t^k.$$
We set $E^k=\cup_{t\in U}E^k_t$ and view it as a (trivial) holomorphic vector bundle over $U$ in a natural way.
If there exists $0<r<s$ such that $B_{r,s}:=\{z\in\mc^m|r\leq\|z\|\leq s\}\subset\Omega_t$ for all $t\in U$ and $\varphi$ is strictly plurisubharmonic on $U\times B_{r,s}$,
then the curvature of the holomorphic Hermitian vector bundle $(E^k,h)$ is strictly positive in the sense of Nakano.
\end{thm}

For Reinhardt domains, we have a similar result which is stronger in form.

\begin{thm}\label{thm:curvartue positive Reinhardt domains from weight}
Let $\Omega\subset U\times\mc^m$ be a family of bounded domains over $U$ that admits a plurisubharmonic defining function,
and $\varphi$ be a $C^2$ plurisubharmonic function defined on some neighborhood of $\overline\Omega$ in $U\times\mc^m$.
We assume that all fibers $\Omega_t\ (t\in U)$ are (connected) Reinhardt domains in $\mc^m$ and $\varphi(t,z)$ is $T^m$- invariant with respect to $z$.
If $\varphi$ is strictly plurisubharmnic on some open subset $O$ in $\Omega$ such that $p(O)=U$,
then for any nonnegative integers $k_1,\cdots, k_m$, the function $\psi(t)$ defined by
$$e^{-\psi(t)}=\int_{\Omega_t}|z_1^{k_1}\cdots z_m^{k_m}|^2e^{-\varphi_t}d\lambda_z$$
is a strictly plurisubharmonic function on $U$.
\end{thm}

We can also deduce some consequences from Theorem \ref{thm:curvartue positive cicular domains from weight} and Theorem \ref{thm:curvartue positive Reinhardt domains from weight} that are parallel to
Corollary \ref{cor:Bergman kernel strict psh} and Corollary \ref{cor:strict minimum principle}.
We leave the details to the readers.

It is helpful to compare Theorem \ref{thm:curvartue positive cicular domains} to Theorem \ref{thm:curvartue positive cicular domains from weight}.
While the strict positivity of the curvature in Theorem \ref{thm:curvartue positive cicular domains} comes from the strict pesudoconvexity of $\Omega$,
it seems that the strict positivity of the curvature in Theorem \ref{thm:curvartue positive cicular domains from weight} essentially comes from the strict plurisubharmonicity of the weight function
on certain subdomain.
In connection to this, we do not know whether Theorem \ref{thm:curvartue positive cicular domains from weight} still holds
if $\Omega$ is just assumed to be a pseudoconvex family (may without plurisubharmonic defining function), with other conditions unchanged.

The second part of the paper aims to establish the connection of Theorem \ref{thm:curvartue positive cicular domains} and Theorem \ref{thm:curvartue positive cicular domains from weight} to the study of ample vector bundles (see \S \ref{sec: ample defn} for definition) over projective manifolds,
which is indeed one of the original motivations for us to consider Theorem \ref{thm:curvartue positive cicular domains} and Theorem \ref{thm:curvartue positive cicular domains from weight}.

By Kodaira's embedding theorem, one can show that a holomorphic vector bundle over a compact complex manifold must be ample
if it is Griffiths positive, i.e., admits a Hermitian metric with positive curvature in the sense of Griffiths.
In 1969, Griffiths conjectured that the converse is true, namely,
such a vector bundle is Griffiths positive if it is ample \cite{Grif69}.
This conjecture is known as \emph{Griffiths conjecture}.
Griffiths conjecture is known to be true if the base space is a Riemannian surface \cite{Ume73}, but is still widely open otherwise.

Along a related direction, Demailly and Skoda proved in 1980 \cite{Dem80} that $E \otimes \det{E} $ is Nakano positive
(i.e., admits a Hermitian metric with positive curvature in the sense of Nakano) if $E$ is ample.
In 2009, Berndtsson proved the following remarkable result.

\begin{thm}[{\cite[Theorem 1.3 and the remark following it]{Ber09}}]\label{thm:ample bunlde positive metric}
If $E$ is an ample vector bundle over a compact complex manifold $X$, then $S^kE \otimes \det{E}$ is Nakano positive for all $k\geq 0$, where $S^kE$
denotes the $k$-th symmetric product of $E$.
\end{thm}

We explain briefly how to deduce Theorem \ref{thm:ample bunlde positive metric} from Theorem \ref{thm:curvartue positive cicular domains} or Theorem \ref{thm:curvartue positive cicular domains from weight}.
Let $E$ be as in Theorem \ref{thm:ample bunlde positive metric}, then one can easily see that there is a smooth strictly negatively curved Finsler metric $h:E^*\ra\mr$ on the
dual bundle $E^*$ of $E$.  Let $\Omega=\{v\in E^*|h(v)<1\}$.
Let $\mb$ be a coordinate ball in $X$ and identify $E^*|_{\mb}$ with $\mb\times\mc^m$ via a local trivialization, where $m$ is the rank of $E$.
Then $\Omega\cap E^*|_{\mb}\subset\mb\times\mc^m$ is a strictly pseudoconvex family of bounded circular domains over $\mb$.
As in Theorem \ref{thm:curvartue positive cicular domains}, we get a Hermitian vector bundle $(E^k, h)$ over $\mb$ whose curvature is strictly positive in the sense of Nakano.
Then the point is that $E^k$ can be canonically identified with $(S^kE\otimes\text{det}E)|_{\mb}$, and the metric $h$ is indeed invariant under transformation of local frame of $E$,
and hence is a global Hermitian metric on $(S^kE\otimes\text{det}E)|_{\mb}$.
In this way, we see that Theorem \ref{thm:ample bunlde positive metric} is a direct consequence of Theorem \ref{thm:curvartue positive cicular domains}.
We can also deduce Theorem \ref{thm:ample bunlde positive metric} from Theorem \ref{thm:curvartue positive cicular domains from weight} in a similar way.
Note that, our argument is different from that in \cite{Ber09} and we do not need to consider the projectivization $\mathbb P(E)$ of $E$.
We hope that our method can throw new light on the study of ample vector bundles and Griffiths conjecture.

We should point out that, to see from a bounded domain $\Omega\subset E^*$ the whole structure of $E$, we have to consider
the natural $S^1$-action on $\Omega$.
So as mentioned above the symmetric structure of $\Omega$ in Theorem \ref{thm:curvartue positive cicular domains} or Theorem \ref{thm:curvartue positive cicular domains from weight}
is indispensable in their application to the proof of Theorem \ref{thm:ample bunlde positive metric}.

We end the introduction by presenting the main ideas of proving Theorem \ref{thm:curvartue positive cicular domains},
Theorem \ref{thm:curvartue positive Reinhardt domains}, and Theorem \ref{thm:curvartue positive cicular domains from weight}.
The two main ingredients in the proofs of the theorems
are Berndtsson's curvature estimate in \cite{Ber09} and Deng-Ning-Wang-Zhou's integral characterization of the Nakano positivity of Hermitian vector bundles \cite{DNWZ}.
More precisely, we first consider product domains and use Berndtsson's estimate to get a positive lower bound of the curvature,
and then take a limit to come back to the metric in the non-product domain case and use Deng-Ning-Wang-Zhou's result to show that
the curvature of the limit metric also has the same lower bound.
The idea in the first step was motivated by Berdtsson's proof of a complex version of the Prekopa's theorem for
non-product domains \cite{Ber98}, and the idea in the second step was applied by Liu-Yang-Zhou to solve a problem
of Lempert via Deng-Ning-Wang-Zhou's result mentioned above \cite{LYZ21}.
In the proofs of Theorem \ref{thm:curvartue positive cicular domains} and Theorem \ref{thm:curvartue positive Reinhardt domains},
one key observation is that a piece of area near the boundary of $\Omega$, no matter how small it is,
can produce a positive lower bound of the curvature of the concerned vector bundle.

Theorem \ref{thm:curvartue positive cicular domains}, Theorem \ref{thm:curvartue positive Reinhardt domains}, Theorem \ref{thm:curvartue positive cicular domains from weight}, and Theorem \ref{thm:curvartue positive Reinhardt domains from weight}
are possible to be generalized to holomorphic vector bundles on more general spaces with general compact group actions.
However, to keep the main idea transparent, we do not touch such general abstract setting.

The remaining of the paper is arranged as follows.
After presenting some necessary preliminaries in \S \ref{sec:preliminary},
we prove Theorem \ref{thm:curvartue positive cicular domains from weight} in \S \ref{sec:Theorem circular domain weight}.
The proof of Theorem \ref{thm:curvartue positive Reinhardt domains from weight} is almost the same as the proof of Theorem \ref{thm:curvartue positive cicular domains from weight}, so we omit it.
We prove Theorem \ref{thm:curvartue positive cicular domains} and Theorem \ref{thm:curvartue positive Reinhardt domains} in \S \ref{sec:thoerem circular and Reinhardt} and deduce the corollaries of them in \S \ref{sec: corollary}.
In the final section \S \ref{sec:ample vector bundle}, we connect Theorem \ref{thm:curvartue positive cicular domains}
and Theorem \ref{thm:curvartue positive cicular domains from weight} with the study of ample vector bundles and deduce from them Theorem \ref{thm:ample bunlde positive metric}.

\subsection*{Acknowledgements}
		The first author thanks Professor Jiafu Ning, Zhiwei Wang, and Xiangyu Zhou for helpful discussions on related topics.
		This research is supported by National Key R\&D Program of China (No. 2021YFA1003100),
        NSFC grants (No. 11871451, 12071310), and the Fundamental Research Funds for the Central Universities.

\section{Preliminaries}\label{sec:preliminary}
In this section, we collect some knowledge that are needed in our discussions.
\subsection{Regular maximum of plurisubharmonic functions}\ \\
Let $\psi \in C^{\infty}(\mr)$ be a nonnegative even function, which is  supported on $[-1, 1]$ and satisfies
$$\int_{\mr} \psi(h) dh=1.$$

\begin{lem}[see {\cite[Lemma (5.18), Chapter I]{Dem}}]\label{lem: Regularization max functions}
	For any $\eta:=(\eta_1, \eta_2) \in (0, +\infty) \times (0, +\infty)$, the function
	$\max\nolimits_{\eta}\colon \mr^2\rw \mr$ defined as
$$(t_1,t_2)\rwo \int_{\mr^2} \max\{t_1+h_1, t_2+h_2\} \frac{1}{\eta_1 \eta_2} \psi(\frac{h_1}{\eta_1}) \psi(\frac{h_2}{\eta_2}) d h_1 d h_2$$
possesses the following properties
\begin{itemize}
\item[(i)] $\max_{\eta}\{t_1, t_2\}$ is non decreasing in all variables, smooth and convex on $\mr^2$;
\item[(ii)] $\max\{t_1, t_2\} \leq  \max_{\eta}\{t_1, t_2\} \leq \max\{t_1 + \eta_1, t_2 + \eta_2\} $;
\item[(iii)] If $u_1, u_2$ are plurisubharmonic functions, then $\max_{\eta}\{u_1, u_2\}$ is also plurisubharmonic.
\end{itemize}

\end{lem}

\subsection{Curvature positivity of Hermitian holomorphic vector bundles}\ \\
Let $X$ be a complex manifold of complex dimension $n$, and $(E, h)$ be a Hermitian holomorphic vector bundle over $X$ of rank $r\leq\infty$.
\par Let $D$ be the $(1,0)$-part of the Chern connection of $(E, h)$, and
\begin{equation}
  \Theta^{(E,h)} := [D, \bar{\partial}] = D \bar{\partial} + \bar{\partial} D
\end{equation}
be the Chern curvature tensor. Over a coordinate chart
$$(\Omega,(t_1,\cdots,t_n))\subset X,$$
we have
$$\partial_{t_j} (u, v) = (D_{t_j} u, v) + (u, \bar{\partial}_{t_j} v),\ \forall u,v \in \Gamma(X,E),$$
where $\partial_{t_j}:= \frac{\partial}{\partial t_j}$ and $\bar{\partial}_{t_j}:= \frac{\partial}{\partial \bar{t}_j}.$
The Chern curvature is
$$\Theta^{(E,h)} =\sum \Theta^{(E,h)}_{jk} dt_j \wedge d\bar{t}_k,$$
where these coefficients are the commutators $\Theta^{(E,h)}_{jk}:=[D_{t_j}, \bar{\partial}_{t_k}]$.\\
\begin{defn}\label{Curvature positivity of Berndtsson}
	The curvature of $(E,h)$ is said to be positive (or strictly positive) in the sense of Nakano if for any nonzero $n$-tuple $(u_1, \cdots  , u_n)$ of sections of $E$
	$$\sum (\Theta_{jk}^{(E,h)} u_j, u_k) \geq 0\ (\text{or} >0).$$
\end{defn}

The following result is obvious.

\begin{lem}[see {\cite[Theorem (14.5), Chapter V]{Dem}}]\label{the direct sum of holomorphic vector bundles}
Let $(F,h)$ be a Hermitian holomorphic vector bundle over $X$, and let $E, G$ be two holomorphic subbundles of $F$ such that
$F=E\oplus G$ and $E$ is orthogonal to $G$, then the curvature of these bundles satisfies
 $$\Theta^F=\Theta^E \oplus \Theta^G.$$
\end{lem}

One of the main ingredients in our argument of the main results is the following result of Berndtsson.

\begin{lem}[{\cite[(3.1)]{Ber09}}]\label{lem: the curvature estimation about vector bundle of holomorphic functions}
If $\Omega$ and $\varphi$ satisfy the conditions in Theorem \ref{thm:Bern direct image bd domain} and $\varphi$ is strictly plurisubharmonic,
then for any smooth sections $u_1, \cdots  , u_n$ of the trivial bundle $E$,
we have
\begin{equation*}
\sum_{j,l} (\Theta_{jl}^E u_j, u_l) \geq \sum_{j,l} \int_D H(\varphi)_{jl} u_j \overline{u_l} e^{-\varphi}d\lambda_z
\end{equation*}
 where
\begin{equation*}
H(\varphi)_{jl}:=\varphi_{jl}- \sum_{\alpha, \beta} \varphi^{\alpha \beta} \varphi_{j \alpha} \overline{\varphi_{l \beta}},
\end{equation*}
where $(\varphi^{\alpha\beta})_{m\times m}$ is the inverse matrix of $(\varphi_{\alpha\beta})_{m\times m}$.
\end{lem}
In the above Lemma, $j,l=1,\cdots, n$ represent the indices of the components of $t=(t_1,\cdots, t_n)$,
$\alpha,\beta=1,\cdots, m$ represent the indices of the components of $z=(z_1,\cdots, z_m)$,
$\varphi_{j,l}=\frac{\partial^2\varphi}{\partial t_j\partial \bar t_l}$ and $\varphi_{j\alpha}, \varphi_{\alpha\beta}$ are given in the same way.

\subsection{Optimal $L^2$-estimate condition and curvature positivity}\

We first recall a fundamental result about the $L^2$-estimate of $\bar\partial$ for a Hermitian holomorphic vector bundle with Nakano positive curvature,
which is due to H\"ormander and Demailly.

\begin{lem}[see {\cite[Theorem (4.5), Chapter VIII]{Dem}}]\label{lem: L2 estimate Nakano}Let $X$ be a complete K\" ahler manifold, with a K\" ahler metric $\omega$ which is not necessarily complete.  Let $(E,h)$ be a Hermitian  vector bundle of rank $r$ over $X$, and assume that the curvature operator $B:=[i\Theta_{E,h},\Lambda_\omega]$ is semi-positive definite everywhere on $\Lambda^{p,q}T_X^*\otimes E$, for some $q\geq 1$. Then for any form $g\in L^2(X,\Lambda^{p,q}T^*_{X}\otimes E)$ satisfying $\bar{\partial}g=0$ and $\int_X\langle B^{-1}g,g\rangle dV_\omega<+\infty$, there exists $f\in L^2(X,\Lambda^{p,q-1}T^*_X\otimes E)$ such that $\bar{\partial}f=g$ and $$\int_X|f|^2dV_\omega\leq \int_X\langle B^{-1}g,g\rangle dV_\omega.$$
\end{lem}

The following result of Deng-Ning-Wang-Zhou shows that the converse of the above Lemma also holds,
and hence gives an equivalent integral form characterization of the curvature positivity of Hermitian holomorphic vector bundles.

\begin{lem}[{\cite[Theorem 1.1]{DNWZ}}]\label{lem: optimal L2-estimate}
Let $U\subset \mc^n$ be a bounded domain, $(E, h)$ be a Hermitian holomorphic vector bundle over $U$ with smooth Hermitian metric $h$,
and $\theta \in C^0 (U, \wedge^{1,1}T^*_U \otimes \operatorname{End}(E))$ with $\theta^*=\theta$.
If for any strictly plurisubharmonic function $\psi$ on $U$ and $f \in C_c^{\infty} (U, \wedge^{n,1}T^*_U \otimes E )$ with $\bar{\partial}f=0$
and $i\partial\bar\partial\psi\otimes Id_E + \theta> 0$ on $\operatorname{supp}(f)$,
there is a measurable section $u$ of $\wedge^{n,0}T^*_U \otimes E $ on $U$, satisfying $\bar{\partial}u=f$ and
\begin{equation}
\int_U |u|_h^2e^{-\psi} d\lambda_z \leq \int_U \langle B_{i\partial\bar\partial\psi, \theta}^{-1}f, f\rangle_{h} d\lambda_z,
\end{equation}
provided that the right hand side is finite, then $i\Theta_{E,h} \geq \theta$ in the sense of Nakano,
 where $\omega=i \sum_{j=1}^n dz_j \wedge d\bar{z}_j$ and
$$B_{i\partial\bar\partial\psi,\theta} = [i\partial\bar\partial\psi\otimes Id_E + \theta, \Lambda_\omega].$$
\end{lem}	

The above Lemma is a modified version of Theorem 1.1 in \cite{DNWZ} (please see {\cite[Remark 1.2]{DNWZ}}.)

\section{The proof of Theorem \ref{thm:curvartue positive cicular domains from weight}}\label{sec:Theorem circular domain weight}
We first give the proof in the case that $\Omega$ is a product domain.

\begin{lem}\label{lem: curvature positivity product domains}
Let $\Omega:=U \times D \subset \mc^n_t \times \mc^m_z$ be a bounded domain, $D$ be a (connected) pseudoconvex circular domain containing the origin.
We assume that $\varphi$ is a $C^2$ plurisubharmonic function defined on some neighborhood of $\overline\Omega$ and is $S^1$-invariant with respect to $z$. Let $k\geq 0$ and $E^k_t$ be the space of homogenous polynomials on $\mc^m$ of degree $k$,  with inner product $h_t$ given by
$$h_t(f,g)=\int_{D}f\bar ge^{-\varphi_t}d\lambda_z,\ f, g\in E_t^k.$$
We set $E^k=\cup_{t\in U}E^k_t$ and view it as a (trivial) holomorphic vector bundle over $U$ in a natural way.
Let $R, M>0$ satisfy
$$\sup\{\|z\|; z\in D\}\leq R,\ \sup\{|\varphi(t,z)|; (t,z)\in\Omega\}\leq M.$$
If there exist $0<r<s$ such that $B_{r,s}:=\{z\in\mc^m|r\leq\|z\|\leq s\}\subset D$ and $\varphi$ is strictly plurisubharmonic on $U\times B_{r,s}$,
then the curvature of the Hermitian holomorphic vector bundle $(E^k,h)$ satisfies:
$$\sum_{j,l} (\Theta_{jl}^{(E^k, h)} u_j, u_l)\geq \delta \sum_j h(u_j,u_j)$$
for any sections $u_1,\cdots, u_n$ of $E^k$,
where $\delta>0$ is a constant depending on $R, M, r, s$ and the complex Hessian of $\varphi$ on $U\times B_{r,s}$.
\end{lem}
\begin{proof}
For any $\epsilon>0,$ let $\varphi_{\epsilon}:=\varphi+ \epsilon (|t|^2+|z|^2)$ and denote the complex Hessian matrix of $\varphi_{\epsilon}$ as $$\left(\begin{array}{clr}
		(\varphi_{\epsilon})_{jl}&(\varphi_{\epsilon})_{j \alpha}\\
		(\varphi_{\epsilon})_{\beta l}&(\varphi_{\epsilon})_{\beta \alpha}
	\end{array}\right),$$
where  $j,l=1,\cdots, n$ represent the indices of the components of $t=(t_1,\cdots, t_n)$,
$\alpha,\beta=1,\cdots, m$ represent the indices of the components of $z=(z_1,\cdots, z_m)$.
Then $\varphi_{\epsilon}$ is strictly plurisubharmonic on $\Omega$.
We consider the Hermitian metric $h^\epsilon$ on $E^k$ given by:
$$h_t^\epsilon(f,g)=\int_{D}f\bar ge^{-(\varphi_\epsilon)_t}d\lambda_z,\ f, g\in E_t^k.$$

Let $E$ be the trivial vector bundle over $U$ as in Theorem \ref{thm:Bern direct image bd domain}.
Then $E^k$ is a holomorphic subbundle of $E$.
Since $D$ is a circular domain containing the origin, any $f\in \mathcal O(D)$ can be represented as  a series
$$f=\sum_{j=0}^{+\infty}f_j$$
that is convergent locally uniformly on $D$, where each $f_j$ is a homogenous polynomial of degree $j$.
For any $S^1$-invariant continuous bounded function $\psi$ on $D$, and any homogenous polynomials $g_j, g_l$ of degree $j$ and $l$ respectively,
we have
$$\int_Dg_j\bar g_l e^{-\psi}=0$$
whenever $j\neq l$.
It follows that, for any $t\in U$, an element $f$ in the orthogonal complement $(E^k)_t^{\bot}$ of $E_t^k$ in $E_t$ has the form
$$f=\sum_{j\geq 0, j\neq k}f_k,$$
where each $f_j$ is a homogeneous polynomial of degree $j$.
Hence $(E^k)_t^{\bot}$ as a vector space is independent of the choice of the weight function $\varphi$ and is also a holomorphic subbundle of $E$.

We now fix an arbitrary $t_0\in U$.
By Lemma \ref{the direct sum of holomorphic vector bundles} and Lemma \ref{lem: the curvature estimation about vector bundle of holomorphic functions},
for any $u_1, \cdots  , u_n$ of $E_{t_0}^k,$ we have
\begin{equation}\label{eqat: curvature estimation Ek}
	\sum_{j,l} (\Theta_{jl}^{(E^k, h^{\epsilon})} u_j, u_l) \geq \int_D\sum_{j,l} H(\varphi_{\epsilon})_{jl}(t_0,z) u_j \overline{u_l} e^{-(\varphi_{\epsilon})_{t_0}}d\lambda_z.
\end{equation}
where $H(\varphi_{\epsilon})$ is a Hermitian matrix defined as in Lemma \ref{lem: the curvature estimation about vector bundle of holomorphic functions}.
\indent

If we write
$$
\left(
\begin{matrix}
  \varphi_{jl}&\varphi_{j \alpha}\\
 \varphi_{\beta l}&\varphi_{\beta \alpha}
\end{matrix}\right)=
\left(
\begin{matrix}
  A & B \\
  C & F
\end{matrix}
\right),
$$
then we have $H(\varphi)=A-BF^{-1}C$ provided that $F$ is nonsigular, and
$$
\left(
\begin{matrix}
  H(\varphi) & 0 \\
  * & F
\end{matrix}
\right)
=\left(
\begin{matrix}
  A-BF^{-1}C & 0 \\
  * & F
\end{matrix}
\right)
=
\left(
\begin{matrix}
  I & -BF^{-1} \\
  0 & I
\end{matrix}
\right)
\left(
\begin{matrix}
  A & B \\
  C & F
\end{matrix}
\right)
\left(
\begin{matrix}
  I & 0 \\
  -(BF^{-1})^* & I
\end{matrix}
\right).
$$
It follows that $H(\varphi)$ is positively definite if $\varphi$ is strictly plurisubharmonic.
So we have
\begin{equation}\label{eqat: curvature estimation via annulus}
\sum_{j,l} (\Theta_{jl}^{(E^k, h^{\epsilon})} u_j, u_l)|_{t_0} \geq \int_{B_{r,s}}\sum_{j,l} H(\varphi_{\epsilon})_{jl}(t_0,z) u_j \overline{u_l} e^{-(\varphi_{\epsilon})_{t_0}}d\lambda_z.
\end{equation}

By assumption and by continuity, there is a constant $\delta_0>0$ such that
$$\sum_{j,l}H(\varphi)_{jl}(t_0,z)u_j\overline{u_l}\geq \delta_0\sum_j |u_j|^2$$
for $z\in B_{r,s}$.
On the other hand, it is clear that $$H(\varphi_{\epsilon})(t_0,z)=H(\varphi)(t_0,z)+o_\epsilon(1)$$
on $B_{r,s}$, where $o_\epsilon(1)$ represents functions on $B_{r,s}$ that converge to 0 uniformly as $\epsilon\ra 0$.
It follows that
$$\sum_{j,l} (\Theta_{jl}^{(E^k, h^{\epsilon})} u_j, u_l)|_{t_0} \geq \delta_0\int_{B_{r,s}}\sum_j(1+o_\epsilon(1))|u_j|^2e^{-(\varphi_\epsilon)_{t_0}}d\lambda_z.$$
Since $h^{\epsilon}$ converges to $h$ in the sense of  $C^2$ as $\epsilon \rw 0^+$, $\Theta_{(E^k, h^{\epsilon})}$ converges to $\Theta_{(E^k, h)}$ as $\epsilon \rw 0^+.$
We thus have
$$\sum_{j,l} (\Theta_{jl}^{(E^k, h)} u_j, u_l)\geq \delta_0\int_{B_{r,s}} \sum_j|u_j|^2e^{-\varphi_{t_0}}d\lambda_z.$$
Note that $u_j$ are homogenous polynomials of degree $k$, $D$ is bounded, and $\varphi(t_0,z)$ is bounded on $\overline D$,
there exists a constant $\delta>0$, which is independent of $u_j$,  such that
$$\delta_0\int_{B_{r,s}} \sum_j|u_j|^2e^{-\varphi_{t_0}}d\lambda_z\geq \delta\int_{D} \sum_j|u_j|^2e^{-\varphi_{t_0}}d\lambda_z.$$
It follows that
$$\sum_{j,l} (\Theta_{jl}^{(E^k, h)} u_j, u_l)\geq \delta\int_D \sum_j|u_j|^2e^{-\varphi_{t_0}}d\lambda.$$
\end{proof}

We shall deduce Theorem \ref{thm:curvartue positive cicular domains from weight} from Lemma \ref{lem: curvature positivity product domains} and Lemma \ref{lem: optimal L2-estimate}.
\begin{thm}[=Theorem \ref{thm:curvartue positive cicular domains from weight}]
	Let $\Omega\subset U\times\mc^m$ be a family of bounded domains over $U$ that admits a plurisubharmonic defining function,
	and $\varphi$ be a $C^2$ plurisubharmonic function defined on some neighborhood of $\overline\Omega$ in $U\times\mc^m$.
	We assume that all fibers $\Omega_t\ (t\in U)$ are (connected) circular domains in $\mc^m$ containing the origin and $\varphi(t,z)$ is $S^1$- invariant with respect to $z$.
	Let $k\geq 0$ and $E^k_t$ be the space of homogenous polynomials on $\mc^m$ of degree $k$, with inner product $h_t$ given by
	$$h_t(f,g)=\int_{\Omega_t}f\bar ge^{-\varphi_t}d\lambda_z,\ f, g\in E_t^k.$$
	We set $E^k=\cup_{t\in U}E^k_t$ and view it as a (trivial) holomorphic vector bundle over $U$ in the natural way.
	If there exist $0<r<s$ such that $B_{r,s}:=\{z\in\mc^m|r\leq\|z\|\leq s\}\subset\Omega_t$ for all $t\in U$ and $\varphi$ is strictly plurisubharmonic on $U\times B_{r,s}$,
	then the curvature of the holomorphic Hermitian vector bundle $(E^k,h)$ is strictly positive in the sense of Nakano.
\end{thm}

\begin{proof}[Proof]
Let $\rho(t,z)$ be a plurisubharmonic defining function of $\Omega$, by averaging, we may assume that $\rho$ is $S^1$-invariant with respect to $z$.
For any fixed $t_0 \in U$ and $0<h<<1$, let $D=\{(t_0, z)\in U\times\mc^n|\rho(t_0,z)\leq h\}$.
Then there exists a neighborhood $U'$ of $t_0$ in $U$ such that $\rho$ and $\varphi$ are defined on some neighborhood of the closure of $U'\times D$ and $p^{-1}(U')\cap\Omega\subset U'\times D$,
where $p:\mc^n\times\mc^m\ra\mc^n$ is the natural projection.
Since the result to be proved is local in nature with respect to $t$, we may assume that $U=U'$, then we have $\Omega\subset\tilde\Omega:=U\times D$.

For any positive integer $N$, let
$$\varphi_N= \varphi +N\max\nolimits_{(\frac{1}{N^2}, \frac{1}{N^2})} \{0, \rho-\frac{1}{N}\}, $$
where $\max\nolimits_{(\frac{1}{N^2}, \frac{1}{N^2})} \{0, \rho\}$ is the regularized max function defined as in Lemma \ref{lem: Regularization max functions}.
For $N>>1$, $\varphi_N$ is equal to $\varphi$ on $\Omega$.
Applying Lemma \ref{lem: curvature positivity product domains} to $\tilde\Omega$ and $\varphi_N$,
we get a constant $\delta>0$ such that
$$\sum (\Theta_{jl}^{(E^k, h^N)} u_j, u_l) \geq \delta \sum \int_D |u_j|^2 e^{-\varphi_{N}}.$$
for any sections $u_1,\cdots, u_n$ of $E^k$,
where the metric $h^N$ on $E^k$ is given by
$$h^N_t(f,g)=\int_{D}f\bar ge^{-(\varphi_N)_t}d\lambda_z,\ f, g\in E_t^k.$$
In other words, if we take
$$\theta=i \delta \sum_j  dt_j \wedge d\bar{t}_j \otimes Id_{E^k} \in C^0 (U, \wedge^{1,1}T^*_U \otimes End(E^k)),$$
then we have $i \Theta_{(E^k, h^N)} \geq_{Nak} \theta$.

We want to apply Lemma \ref{lem: optimal L2-estimate} to prove that $i \Theta_{(E^k, h)} \geq_{Nak} \theta$.
The main idea is as follows.
From the above curvature estimate and the $L^2$-estimate of $\bar\partial$,
we know that $(E^k, h^N)$ satisfy the $L^2$-estimate condition presented in Lemma \ref{lem: optimal L2-estimate}.
As $N\ra\infty$, we have $h^N\ra h$ and one can see that $(E,h)$ also satisfies the $L^2$-estimate condition.
Then it follows from Lemma \ref{lem: optimal L2-estimate} that the curvature of $(E,h)$ satisfies $i \Theta_{(E^k, h)} \geq_{Nak} \theta$.
The detail of the argument is as follows.

Let $\psi(t)$ be a strictly plurisubharmonic function on $U$,
and $f \in C_c^{\infty} (U, \wedge^{n,1}T^*_U \otimes E^k)$ satisfies $\bar{\partial}f=0$ and
$$\int_U <B_{i\partial\bar\partial\psi, \theta}^{-1}f, f>_{h}e^{-\psi} d\lambda_t < + \infty,$$
where $\omega=i \sum_{j=1}^n dt_j \wedge d\bar{t}_j$ and $B_{i\partial\bar\partial\psi, \theta}$ is given as in Lemma \ref{lem: optimal L2-estimate}.
Then there exists $M>0$ such that
$$\int_U <B_{i\partial\bar\partial\psi, \theta}^{-1} f, f>_{h^N}e^{-\psi} d\lambda_t \leq M,\ \forall N.$$

By Lemma \ref{lem: L2 estimate Nakano},
there are measurable sections $u_N$ of  $\wedge^{n,0}T^*_U \otimes E^k$ on $U$, such that $\bar{\partial}u_N=f$ and
$$\int_U |u_N|^2_{h^N} e^{-\psi} d\lambda_t \leq \int_U <B_{i\partial\bar\partial\psi, \theta}^{-1}f, f>_{h^N} e^{-\psi} d\lambda_t \leq M.$$
Since $\varphi^N$ and $\varphi$ are equal on $\Omega$, we have
$$\int_U |u_N|^2_{h} e^{-\psi} d\lambda_t\leq \int_U |u_N|^2_{h^N} e^{-\psi} d\lambda_t\leq M$$
for all $N\geq 1$.
In particular, $\{u_N\}$ is a bounded sequence in the Hilbert space $H$ of square integrable sections of $\wedge^{n,0}T^*_U \otimes E^k$ on $U$ with weight $e^{-\psi}$. Hence there is a subsequence of $\{u_N\}$, assumed to be $\{u_N\}$ itself without loss of generality, that converges weakly in $H$ to some $u$.
Note that we also have $\bar\partial u=f$ in the sense of distribution.
On one hand, we have
$$\int_U|u|^2_he^{-\psi} d\lambda_t\leq\limsup_{N\ra\infty}\int_U |u_N|^2_{h} e^{-\psi} d\lambda_t,$$
and on the other hand,
we have
$$\lim_{N\ra\infty}\int_U <B_{i\partial\bar\partial\psi, \theta}^{-1}f, f>_{h^N} e^{-\psi} d\lambda_t=
\int_U <B_{i\partial\bar\partial\psi, \theta}^{-1}f, f>_{h} e^{-\psi} d\lambda_t$$
by Lebesgue's dominated convergence theorem.
So we get
$$\int_U|u|^2_he^{-\psi} d\lambda_t\leq \int_U <B_{i\partial\bar\partial\psi, \theta}^{-1}f, f>_{h} e^{-\psi} d\lambda_t.$$
It follows from Lemma \ref{lem: optimal L2-estimate} that $i \Theta_{(E^k, h)} \geq_{Nak} \theta$.
\end{proof}

\section{The proof of Theorem \ref{thm:curvartue positive cicular domains} and Theorem \ref{thm:curvartue positive Reinhardt domains}}\label{sec:thoerem circular and Reinhardt}
The difficulty of Theorem \ref{thm:curvartue positive cicular domains} compared with Theorem \ref{thm:curvartue positive cicular domains from weight}
is that the weight function does not have strict plurisubharmonicity.
We will use the strict psedoconvexity of the domain to get the Nakano positivity.
In the proof of Theorem \ref{thm:curvartue positive cicular domains},
in addition to using Berndtsson's estimate of curvature (Lemma \ref{lem: the curvature estimation about vector bundle of holomorphic functions}) and Deng-Ning-Wang-Zhou's integral characterization of the Nakano positivity of Hermitian vector bundles (Lemma \ref{lem: optimal L2-estimate}),
an important role is also played by the simple observation that the integral $\int^r_0Ne^{-Nh(x)}dx$ has a uniform positive limit as $N\ra\infty$ for all $r>0$ and all smooth function $h$ with $h(0)=0$ and $h'(0)\leq 1$.

We first give a Lemma.
\begin{lem}\label{lem:volume form}
Let $\Omega$ be a bounded domain in $\mr^n$ with $C^2$- boundary.
For any $0<r<<1$, let $$\Omega_r:=\{x\in\mr^n\backslash\Omega|d(x,\partial\Omega)<r\}.$$
Then there exists a constant $c>0$ such that
$$\int_{\Omega_r} h dx_1 \wedge \cdots \wedge dx_n \geq c \int_{\partial \Omega} dS \int_{0}^r h(\zeta+t\mathbf{n}_\zeta) dt$$
for any positive integrable functions $h$ on $\Omega_r$,
where $\mathbf{n}_\zeta$ is the outward unit normal of $\partial\Omega$ at $\zeta$ and $dS$ is the volume form on $\partial \Omega$.
\end{lem}
\begin{proof}
  We can choose $r_0>0$ such that the map
  $$f:\partial\Omega\times[0,r_0)\ra\Omega_{r_0}; (\zeta,t)\mapsto \zeta+t\mathbf{n}_\zeta$$
  is a diffeomorphism.
  Let $\mu=dS\wedge dt$ be the product measure on $\partial\Omega\times[0,r_0)\ra\Omega_{r_0}$ and $\mu_0$ be the Lebesgue measure on $\Omega_{r_0}$.
  Then there is a continuous positive function $\sigma$ on $\Omega_{r_0}$ such that $\mu_0=\sigma\cdot f_*\mu$ on $\Omega_{r_0}$.
  For any $0<r<r_0$, taking $c=\min\{\sigma(x)|x\in\Omega_{r}\}$, then $c>0$ and $\mu_0\geq c f_*\mu$.
 From it the lemma follows.
\end{proof}

\begin{thm}[=Theorem \ref{thm:curvartue positive cicular domains}]
Let $\Omega\subset U\times\mc^m$ be a strictly pseudoconvex family of bounded domains over $U\subset\mc^n$ and $\varphi$ be a $C^2$ plurisubharmonic function defined on some neighborhood of $\overline\Omega$ in $U\times\mc^m$.
We assume that all fibers $\Omega_t\ (t\in U)$ are (connected) circular domains in $\mc^m$ containing the origin and $\varphi(t,z)$ is $S^1$- invariant with respect to $z$.
Let $k\geq 0$ and $E^k_t$ be the space of homogenous polynomials on $\mc^m$ of degree $k$,  with inner product $h_t$ given by $$h_t(f,g)=\int_{\Omega_t}f\bar ge^{-\varphi_t}d\lambda_z,\ f, g\in E_t^k.$$
We set $E^k=\cup_{t\in U}E^k_t$ and view it as a (trivial) holomorphic vector bundle over $U$ in a natural way.
Then the curvature of the holomorphic Hermitian vector bundle $(E^k,h)$ is strictly positive in the sense of Nakano.
\end{thm}

\begin{proof}
Since $\Omega$ is strictly pseudoconvex with $C^2$ boundary, there is a defining function $\rho$
that is strictly plurisubharmonic on some neighborhood $\tilde\Omega$ of $\overline{\Omega}$ in $U\times\mc^m$ and $S^1$ invariant with respect to $z$.

For any fixed $t_0 \in U$, we can take a neighborhood $U'$ of $t_0$ in $U$ and a pseudoconvex circular domain $D\subset\mc^m$
such that $p^{-1}(U')\cap\Omega\subset U'\times \overline D\subset\tilde\Omega$.
Since the conclusion to be proved is local in nature on $t$,
we may assume that $U=U'$.
We denote $U\times D$ by $\Omega'$.

For $ N \in \mathbb{Z}_+$, we set
$$\varphi_N= \varphi +N\max\nolimits_{(\frac{1}{N^3}, \frac{1}{N^3})} \{0, \rho\}, $$
which is a $C^2$ plurisubharmonic function
defined on $\tilde\Omega$ and is $S^1$- invariant with respect to $z$, where $\max_{(\frac{1}{N^3}, \frac{1}{N^3})} \{0, \rho\}$ is the regularized max function
defined as in Lemma \ref{lem: Regularization max functions}.
For any $\epsilon>0,$ define
   $$\varphi_{N, \epsilon}:=\varphi+N\max\nolimits_{(\frac{1}{N^3},\frac{1}{N^3})}\{0,\rho\}
   +\epsilon|t|^2+\epsilon|z|^2.$$
Let $h^{N,\epsilon}$ be the Hermitian metric on $E^k$ given by
$$h^{N,\epsilon}_t=\int_Df\bar ge^{-(\varphi_{N,\epsilon})_t}d\lambda_z,\ f,g\in E_t.$$

By Lemma \ref{lem: the curvature estimation about vector bundle of holomorphic functions},
we know for any $u_1, \cdots  , u_n$ of $E_{t_0}^k$ that
$$\sum_{j,l} (\Theta_{jl}^{(E^k, h^{N,\epsilon})} u_j, u_l) \geq \int_D \sum_{j,l} H(\varphi_{N,\epsilon})_{jl}(t_0,z) u_j \overline{u_l} e^{-(\varphi_{N,\epsilon})_{t_0}}d\lambda,$$
where $H(\varphi_{N,\epsilon})_{jl}$ is defined as in Lemma \ref{lem: the curvature estimation about vector bundle of holomorphic functions}.

For $0<r<<1$, as in Lemma \ref{lem:volume form}, we set
$$\Omega_{t_0,r}=\{z\in\mc^m\backslash\Omega_{t_0}|d(z,\partial\Omega_{t_0})<r\}$$
and set $\Omega^N_{t_0,r}=\Omega_{t_0,r}\backslash \Omega_{t_0,1/N^2}$ for $N>0$.
We now fix such an $r$ such that $\Omega_{t_0,r}\subset D$.
Note that $\max\nolimits_{(\frac{1}{N^3},\frac{1}{N^3})}\{0,\rho\}=\rho$ on $\Omega^N_{t_0,r}$ for all $N$.

Note that
$$H(\varphi+N\rho +\epsilon|t|^2+\epsilon|z|^2)=NH(\rho+(\varphi+\epsilon|t|^2+\epsilon|z|^2)/N),$$
we have
$$H(\varphi+N\rho +\epsilon|t|^2+\epsilon|z|^2) \geq \frac{N}{2} H(\rho)$$
on $\Omega^N_{t_0,r}$ for $N$ sufficiently large.
Combining with Lemma \ref{lem:volume form}, we can see there exist constants $\delta_0, \delta_1>0$ such that
\begin{align*}
	&\sum_{j,l} (\Theta_{jl}^{(E^k, h^{N, \epsilon})} u_j, u_l)\\
    \geq &\int_D \sum H(\varphi_{N,\epsilon})_{jl}(t_0,z) u_j \overline{u_l} e^{-(\varphi_{N,\epsilon})_{t_0}}d\lambda_z\\
	\geq &\int_{\Omega^N_{t_0,r}} \sum H(\varphi_{N,\epsilon})_{jl}(t_0,z) u_j \overline{u_l} e^{-(\varphi_{N,\epsilon})_{t_0}}d\lambda_z\\
	\geq & \delta_0\int_{\Omega^N_{t_0,r}}N\sum |u_j|^2 e^{-N\rho}d\lambda_z\\
    \geq & \delta_1 \int_{\zeta\in\partial \Omega_{t_0}} dS\int^r_{1/N^2} \sum N |u_j(\zeta+\tau \mathbf{n}_\zeta)|^2 e^{-N\rho(\zeta+\tau \mathbf{n}_\zeta)}d\tau\\
     \geq & \delta_1 \int_{\zeta\in\partial \Omega_{t_0}} dS \sum \inf_{1/N^2\leq \tau\leq r} |u_j(\zeta+\tau \mathbf{n}_\zeta)|^2 \int^r_{1/N^2} \sum N  e^{-N\rho(\zeta+\tau \mathbf{n}_\zeta)}d\tau\\
      \geq & \delta_1 \int_{\zeta\in\partial \Omega_{t_0}} dS \sum\inf_{0\leq \tau\leq r} |u_j(\zeta+\tau \mathbf{n}_\zeta)|^2 \int^r_{1/N^2} N  e^{-NT\tau}d\tau,
\end{align*}
where $\mathbf{n}_\zeta$ is the unit outward normal of $\partial\Omega_{t_0}$ at $\zeta$,
$dS$ is the volume form on $\partial\Omega_{t_0}$,
and $T>0$ is a constant such that $\rho(\zeta+\tau \mathbf{n}_\zeta)\leq T\tau$ for all $\zeta\in\partial\Omega_{t_0}$ and $0\leq\tau\leq r$.
We now need the obvious but important fact that $\lim_{N\ra\infty}\int^r_{1/N^2} N  e^{-NT\tau}d\tau=\frac{1}{T}>0$.
We then get from the above calculation that
\begin{align*}
	\sum_{j,l} (\Theta_{jl}^{(E^k, h^{N, \epsilon})} u_j, u_l)
\geq \delta_2 \sum\int_{\partial\Omega_{t_0}} \inf_{0\leq \tau\leq r} |u_j(\zeta+\tau \mathbf{n}_\zeta)|^2 dS
\end{align*}
for some constant $\delta_2>0$ and for $N$ sufficiently large.
Let $\epsilon\ra 0$, and denote $h^{N,0}$ by $h^N$, we get
\begin{align}\label{eq:curvature estimate vir boundary}
	\sum_{j,l} (\Theta_{jl}^{(E^k, h^{N})} u_j, u_l)
\geq \delta_2 \sum\int_{\partial\Omega_{t_0}} \inf_{0\leq \tau\leq r} |u_j(\zeta+\tau \mathbf{n}_\zeta)|^2 dS
\end{align}
for $N$ sufficiently large.

For $u\in E^k_{t_0}$, we need to control its norm
$$\|u\|_{h^N_{t_0}}=\int_D|u|^2e^{-(\varphi_N)_{t_0}}d\lambda_z$$
in terms of the integral $\int_{\partial\Omega_{t_0}} \inf_{0\leq \tau\leq r} |u(\zeta+\tau \mathbf{n}_\zeta)|^2 dS$,
where $\varphi_N=\varphi_{N,0}$.

Let $Q=\{u\in E^k_{t_0}; \|u\|^2_{h^N}= 1\}$.
Note that functions in $Q$ are homogenous polynomials of degree $k$ and $\Omega_{t_0}$ contains the origin,
we can choose a constant $M>0$ and a large ball $B$ with $\overline D\subset B$ such that
$\int_B|u|^2d\lambda_z\leq M$ for all $u\in Q$.
By Cauchy's inequality for holomorphic functions,
there is a constant $C>0$ such that
$|du^2|<C$ on $D$ for all $u\in Q$.
It follows that
\begin{equation}\label{eq:curvature estimate via b}
\inf_{0\leq \tau\leq r} |u(\zeta+\tau \mathbf{n}_\zeta)|^2\geq |u(\zeta)|^2-rC
\end{equation}
for all $\zeta\in\partial\Omega_{t_0}$ and for all $u\in Q$.

We now move to prove that we can choose $r$ and a constant $\delta_3 >0$ such that
$$\int_{\partial \Omega_{t_0}} \inf_{0\leq \tau\leq r} |u(\zeta+\tau \mathbf{n}_\zeta)|^2 dS \geq \delta_3$$
for all $u\in Q$.
By the maximum principle and continuity, we can take $\zeta'\in \partial\Omega_{t_0}$ such that
$|u|$ takes its maximum on $\overline\Omega_{t_0}$ at $\zeta'$.
Again, since functions in $Q$ are homogenous polynomials of degree $k$ and $\Omega_{t_0}$ contains the origin,
we can choose a constant $C_1>0$ such that
$\int_{\Omega_{t_0}}|u|^2d\lambda_z\geq C_1$ for all $u\in Q$.
It follows that
$$|u(\zeta')|^2 \geq \frac{C_1}{|\Omega_{t_0}|},$$
where $|\Omega_{t_0}|$ is the Lebesgue measure of $\Omega_{t_0}$.

Again by Cauchy's inequality, if choosing $0<r<\frac{C_1}{2C|\Omega_{t_0}|}$, we get
$$|u(\zeta)|^2 \geq |u(\zeta')|^2-Cr \geq \frac{C_1}{2|\Omega_{t_0}|} $$
for all $u\in Q$ and for all $\zeta\in\partial\Omega_{t_0}$ with $|\zeta-\zeta'|<r$.
It follows that
\begin{equation*}
\begin{split}
&\int_{\partial\Omega_{t_0}}\inf_{0\leq \tau\leq r} |u(\zeta+\tau \mathbf{n}_\zeta)|^2dS\\
\geq & \int_{B(\zeta',r) \cap \partial \Omega_{t_0}}\inf_{0\leq \tau\leq r} |u(\zeta+\tau \mathbf{n}_\zeta)|^2dS\\
\geq &  \int_{B(\zeta',r) \cap \partial \Omega_{t_0}} (|u|^2-rC)dS\\
\geq & \frac{C_1}{2|\Omega_{t_0}|} |B(\zeta',r) \cap \partial \Omega_{t_0}| ,
\end{split}
\end{equation*}
where $B(\zeta',r)$ is the ball in $\mc^m$ with center $\zeta'$ and radius $r$.
Note that $\partial\Omega_{t_0}$ is compact and the function
$$\sigma : \partial\Omega_{t_0} \longrightarrow \mr : \zeta \rightarrow |B(\zeta,r) \cap \partial\Omega_{t_0}|$$
is continuous and positive, we have
 $$\delta_3:=\inf_{\zeta \in \partial D_r} |B(\zeta,r) \cap \partial\Omega_{t_0}|>0.$$
 So we get
 $$\int_{\partial \Omega_{t_0}} \inf_{0\leq \tau\leq r} |u(\zeta+\tau \mathbf{n}_\zeta)|^2 dS \geq \delta_3.$$
By \eqref{eq:curvature estimate vir boundary}, for $N$ sufficiently large, we have
\be\label{eq:curvature of hN}
\sum_{j,l} (\Theta_{jl}^{(E^k, h^{N})} u_j, u_l)\geq\delta\sum_j\|u_j\|^2_{h_{t_0}^N}
\ee
for any tuple $u_1,\cdots, u_n\in E^k_{t_0}$.
Just as the last step in the proof of Theorem \ref{thm:curvartue positive cicular domains from weight},
we can derive from \eqref{eq:curvature of hN} and Lemma \ref{lem: optimal L2-estimate} that
$$\sum_{j,l} (\Theta_{jl}^{(E^k, h)} u_j, u_l)\geq\delta\sum_j\|u_j\|^2_{h_{t_0}}$$
for any tuple $u_1,\cdots, u_n\in E^k_{t_0}$.
In particular, the curvature of $(E^k,h)$ is strictly positive in the sense of Nakano.
\end{proof}	

Similar results holds for a strictly pseudoconvex family of Reinhardt domains.

\begin{thm}[=Theorem \ref{thm:curvartue positive Reinhardt domains}]
	Let $\Omega\subset U\times\mc^m$ be a strictly pseudoconvex family of bounded domains over $U\subset\mc^n$ and $\varphi$ be a $C^2$ plurisubharmonic function defined on some neighborhood of $\overline\Omega$ in $U\times\mc^m$.
	We assume that all fibers $\Omega_t\ (t\in U)$ are (connected) Reinhardt domains in $\mc^m$ and $\varphi(t,z)$ is $T^m$ invariant with respect to $z$.
	Then for any nonnegative integers $k_1,\cdots, k_m$, the function $\psi(t)$ defined by
	$$e^{-\psi(t)}=\int_{\Omega_t}|z_1^{k_1}\cdots z_m^{k_m}|^2e^{-\varphi_t}d\lambda_z$$
	is a strictly plurisubharmonic function on $U$.
\end{thm}

\begin{proof}
Since the proof is almost the same as the proof of Theorem \ref{thm:curvartue positive cicular domains}, we just give a sketch of it.

For any nonnegative integers $k_1,\cdots, k_m$, we consider the 1-dimensional vector space
$$E^{k_1,\cdots, k_m}_t =\mc z_1^{k_1}\cdots z_m^{k_m},$$
with inner product $h_t$ given by
	$$h_t(f,g)=\int_{\Omega_t}f\bar ge^{-\varphi_t}d\lambda_z,\ f, g\in E^{k_1,\cdots, k_m}_t.$$
	We set $E^{k_1,\cdots, k_m}=\cup_{t\in U}E^{k_1,\cdots, k_m}_t$ and view it as a holomorphic line bundle over $U$ in the natural way.

Since the conclusion to be proved is local in nature with respect to $t\in U$,
we may assume there is a bounded pseudoconvex Reinhardt domain $D\subset\mc^n$ such that $\Omega\subset \Omega':=U\times D$
and $\varphi$ and $\rho$ are defined on some neighborhood of $\overline{\Omega'}$.

Note that
$$\int_D z_1^{k_1}\cdots z_m^{k_m}\overline{z_1^{l_1}\cdots z_m^{l_m}}e^{-\varphi_t}d\lambda_z=0$$
for any nonnegative integers $k_1,\cdots, k_m$ and $l_1,\cdots, l_m$ with $k_j\neq l_j$ for some $1\leq j\leq m$.
So by Lemma \ref{the direct sum of holomorphic vector bundles} the curvature of $(E^{k_1,\cdots, k_m},h)$ is the restriction of the curvature of
$(E, h')$ on $E^{k_1,\cdots, k_m}$,
where $(E, h')$ represents the vector bundle given in Theorem \ref{thm:Bern direct image bd domain}
with $\Omega$ replaced by $\Omega'$.

With the above discussions at hand,
the remaining  of the proof of the theorem can go ahead following the same way as in the proof of Theorem \ref{thm:curvartue positive cicular domains},
and we omit the details here.
\end{proof}


\section{Some consequences of Theorem \ref{thm:curvartue positive cicular domains} and Theorem \ref{thm:curvartue positive Reinhardt domains}}\label{sec: corollary}
We now discuss some consequences of Theorem \ref{thm:curvartue positive cicular domains} and Theorem \ref{thm:curvartue positive Reinhardt domains}.
\subsection{Consequences in complex analysis}\
\par
We prove Corollary \ref{cor:Bergman kernel strict psh} and Corollary \ref{cor:strict minimum principle} in this subsection.
\begin{cor}[=Corollary \ref{cor:Bergman kernel strict psh}]
	Let $\Omega\subset U\times\mc^m$ be a strictly pseudoconvex family of bounded domains over $U\subset\mc^n$ and $\varphi$ be a $C^2$ plurisubharmonic function defined on some neighborhood of $\overline\Omega$ in $U\times\mc^m$
	that satisfy the conditions in  Theorem \ref{thm:curvartue positive cicular domains} or Theorem \ref{thm:curvartue positive Reinhardt domains}.
	For $t\in U$, let $K(t,z)$ be the weighted Bergman kernel of $\Omega_t$ with weight $\varphi_t$.
	Then $\ln K(t,z)$ is a strictly plurisubharmonic function on $\Omega$.	
\end{cor}

The proof is provided in the following discussion, which indeed gives us more information.

We assume $\Omega$ and $\varphi$ satisfies the conditions in Theorem \ref{thm:curvartue positive cicular domains},
and the remaining case can be proved in the same way.

Note that $\ln K(t,z)$ is strictly plurisubharmonic with respect to $z$,
it is enough to prove that for any $(t_0,z_0)\in \Omega$ and any local holomorphic map
$\xi(t):B\ra\mc^m$ defined on some small neighborhood $B$ of $t_0$ with $\xi(t_0)=z_0$,
the function $\ln K(t, \xi(t))$ is strictly plurisubharonic as a function on $B$
(the reason is that any non-vertical tangent vector of $\Omega$ at $(t_0, z_0)$ lies in the image
of $d\xi(t_0))$ for some such a map $\xi$).

Let $E^k_t$ be the space with inner product defined as in Theorem \ref{thm:curvartue positive cicular domains}, and let $u^k_1,\cdots, u^k_{m_k}$ be an orthogonal
normal basis of $E^k_t$.
We set
$$K^k(t,z)=\sum^{m_k}_{j=1}|u_j(z)|^2,$$
then it is clear that
\be\label{eq:Bergman kernel decomposition}
K(t,z)=\sum^\infty_{k=0}K^k(t,z).
\ee

Let $p:\Omega\ra U$ be the natural projection.
Then the pull back
$$(\tilde E^k,\tilde h):=(p^*E^k, p^*h)$$
of the bundle $(E^k,h)$ on $U$
is a Hermitian holomorphic vector bundle over $\Omega$ whose curvature is semi-positive in the sense of Nakano.

Let $F=\Omega\times\mc$ be the trivial line bundle on $\Omega$ and denote by $e$ the canonical frame of $F$ on $\Omega$.
Then we have a natural vector bundle morphism
$\sigma_k: \tilde E^k\ra L$ given by
$$f\mapsto (t,z, f(z))\in F$$
for $f\in\tilde E^k_{(t,z)}=E^k_t$.
Let $$\Omega^k=\{(t,z)\in \Omega| K^k(t,z)\neq 0\},$$
or equivalently, $(t,z)\in \Omega^k$
if and only if $f(z)\neq 0$ for some homogenous polynomial $f$ on $\mc^m$ of degree $k$.
Then $\sigma_k$ is a surjective bundle morphism from $\tilde E^k|_{\Omega^k}$ to $F|_{\Omega^k}$.
One can see that the quotient metric, say $h^k$ on $F|_{\Omega^k}$ induced from this morphism is given by
$$\|e\|^2_{h^k}=\frac{1}{K^k(t,z)}=e^{-\ln K^k(t,z)}.$$
Since the curvature of $(\tilde E^k, \tilde h)$ is semi-positive in the sense of Nakano,
and note the curvature increasing property under taking quotient metric {\cite[a) in Proposition (6.10)]{Dem}},
we know the curvature of $(F|_{\Omega^k}, h^k)$ is semi-positive,
which implies that $\ln K^k(t,z)$ is  plurisubharmonic on $\Omega$.

For any given $(t_0,z_0)\in \Omega$, and any holomorphic map
$\xi(t):B\ra\mc^m$ defined on some small neighborhood $B$ of $t_0$ with $\xi(t_0)=z_0$,
we denote by
$$\Gamma=\{(t,\xi(t))|t\in B\}\subset\Omega$$
the graph of $\xi$.
Then $(\tilde E^0,\tilde h)|_{\Gamma}$ is a (trivial) Hermitian line bundle over $\Gamma$ whose curvature is strictly positive,
since $p|_\Gamma:\Gamma\ra B$ is a biholomoprhic map.
Note also that $\sigma_0:\tilde E^0\ra L$ is an isomorphism of vector bundles,
it follows that $\ln K^0(t, \xi(t))$ is strictly plurisubharmonic on $\Gamma$,
and hence is strictly plurisubharmonic as a function of $t$.
By \eqref{eq:Bergman kernel decomposition}, we know that $\ln K(t,\xi(t))$ is strictly plurisubharmonic as a function of $t$.
Hence $\ln K(t,z)$  is strictly plurisubharmonic on $\Omega$.
The proof of the above corollary is complete.

In fact, by the same argument, one can show, for any nonnegative integer $k$, that "the relative log character Bergman kernel" $\ln K^k(t,z)$
is plurisubharmonic on $\Omega$ and is strictly plurisubharmonic on $\Omega^k$.

\begin{cor}[=Corollary \ref{cor:strict minimum principle}]
	Let $\Omega\subset U\times\mc^m$ be a strictly pseudoconvex family of domains over $U\subset\mc^n$ and $\varphi$ be a $C^2$ plurisubharmonic function defined on some neighborhood of $\overline\Omega$ in $U\times\mc^m$.
	\bi
	\item[(1)] If $\Omega$ and $\varphi$ satisfy the conditions in  Theorem \ref{thm:curvartue positive cicular domains} or Theorem \ref{thm:curvartue positive Reinhardt domains}, then the function $\tilde\varphi$ defined by
	$$e^{-\tilde\varphi(t)}=\int_{\Omega_t}e^{-\varphi(t,z)}d\lambda_z$$
	is a strictly plurisubharmonic function on $U$.
	\item[(2)]If all fibers $\Omega_t$ are tube domains of the form $X_t+i\mr^m$ with$X_t$ bounded, and $\varphi(t,z)$ does not depend on the imaginary part of $z$,
	then the function $\tilde\varphi$ defined by
	$$e^{-\tilde\varphi(t)}=\int_{X_t}e^{-\varphi(t,\text{Re}z)}d\lambda_{\text{Re}z}$$
	is a strictly plurisubharmonic function on $U$.
	\ei
\end{cor}
\begin{proof}
It is clear that (1) is equivalent to the curvature strict positivity of $(E^0,h)$ in Theorem \ref{thm:curvartue positive cicular domains} or Theorem \ref{thm:curvartue positive Reinhardt domains}. We now give the proof of (2).

Let us consider the map
	\begin{equation*}
		\begin{split}
			f: \Omega &\rightarrow \mc^n_t\times\mc^m_w\\
			(t_1,\cdots,t_n, z_1,\cdots, z_m)& \mapsto (t_1,\cdots, t_n, e^{z_1}, \cdots, e^{z_m}),
		\end{split}
	\end{equation*}
then $\Omega^{*}:=f(\Omega)\subset\mc^n_t\times\mc^m_w$ is a strictly pseudoconvex family of Reinhardt domains over  $U$.
Note that
$$\psi(t,w) := \varphi(t, \ln|w_1|,\cdots, \ln|w_m|)+2(\ln|w_1|+\cdots+\ln|w_m|)$$
is a $C^2$ and plurisubharmonic function defined on some neighborhood of the closure of $\Omega^*$ in $U\times\mc^m$,
applying (1) to $\Omega^*$ and $\psi$, we see that the function $\tilde\varphi$ defined by
	$$e^{-\tilde\varphi(t)}=\int_{X_t}e^{-\varphi(t,\text{Re}z)}d\lambda_{\text{Re}z}=\frac{1}{(2\pi)^n}\int_{\Omega^*_t}e^{-\psi(t,w)}d\lambda_w $$
	is a strictly plurisubharmonic function on $U$.
\end{proof}

\subsection{Consequences in convex analysis}\
\par
The bridge connecting strictly convex families of bounded domains in $\mr^m$ and strictly pseudoconvex families of tube domains in $\mc^m$
is indicated in the proof of the following corollary.

\begin{cor}[=Corollary \ref{cor:strict prekopa}]
Let $D\subset U_0\times\mr^m$ be a strictly convex family of bounded domains over a domain $U_0\subset\mr^n$ and $\varphi$ be a $C^2$ convex function defined on some neighborhood of the closure of $D$ in $U_0\times\mr^m$.
Then the function $\tilde\varphi$ defined by
$$e^{-\tilde\varphi(t)}=\int_{D_t}e^{-\varphi(t,x)}d\lambda_x$$
is a strictly convex function on $U_0$.
\end{cor}
\begin{proof}
We first complexify $U_0$ to $U=U_0\times i\mr^n_{l}$ with complex coordinate $\tau=t+il$,
then $U$ is a domain in $\mc^n_{\tau}$.
We secondly complexify $\mr^m$ to $\mr_x^m+i\mr_y^m=\mc^m_z$, with complex coordinate $z=x+iy$.
Then
$$\Omega=D+i\mr^{n+m}=\{(\tau, z)\in \mc^n\times\mc^m|(\text{Re}\tau, \text{Re}z)\in D\}$$
is a strictly pseudoconvex family of tube domains over $U$.
For $\tau\in U$, $\Omega_{\tau}$ is a tube domain of the form $\Omega_\tau=D_\tau+i\mr^m$,
where $D_\tau\subset\mr^m$ can be naturally identified with $D_{\text{Re}\tau}$.
By setting
$$\psi(\tau,z)=\varphi(\text{Re}\tau, \text{Re}z),$$
we extend $\varphi$ to a $C^2$ plurisubharmonic function $\psi$ on some neighborhood of $\overline\Omega$ in $U\times\mc^m$,
such that $\psi(\tau,z)$ is independent of the imaginary part of $\tau, z$.
By (2) in Corollary \ref{cor:strict minimum principle}, the function $\tilde{\psi}$ defined by
$$e^{-\tilde\psi(\tau)}=\int_{D_\tau}e^{-\psi(\tau,\text{Re}z)} d\lambda_{\text{Re}z}$$
is a strictly plurisubharmonic function on $U$.
	It is clear that $\tilde{\psi}(\tau)$ is independent of the imaginary part of $\tau$ and $\tilde{\psi}|_{U_0}=\tilde{\varphi}$,
	thus $\tilde\varphi$ is a strictly convex function on $U_0$.
\end{proof}	

\subsection{Curvature negativity of determinant line bundle}\
\par	
We now explain the meaning of Corollary \ref{cor: detE strict negative} and give its proof.

Let $\pi:E\ra X$ be a holomorphic vector bundle of rank $m$ over a complex manifold $X$ equipped with a smooth Finsler metric $h$.
By definition,  $h$ is a continuous function $h:E\ra \mr$ such that $h\geq 0$,
$h(\lambda v)=|\lambda|h(v)$ for $\lambda\in\mc$ and $v\in E$, and $h$ is smooth on $E\backslash Z_E$, where $Z_E\subset E$ is the zero section of $E$.
Recall that $(E,h)$ is defined to be strictly negatively curved if $\ln h$ is strictly plurisubharmonic on $E\backslash Z_E$.

We now define the Hermitian metric $\text{det}h$ induced from $h$ on the determinant line bundle $\text{det}E=\Lambda^mE$ of $E$
via the measure $\mu$ on $E_t$ with $\mu(B_t)=1$ for $t\in X$,
where $$B_t=\{v\in E_t|h(v)\leq 1\}.$$
A more explicit description of $\det h$ in terms of local frame is as follows.
Let $e_1,\cdots, e_m$ be a holomorphic local frame of $E$ over some open set $U\subset X$.
We get a local trivialization of $E$ over $U$:
$$\phi:E|_U\ra U\times\mc^m, \ (t,z_1v_1+\cdots+z_rv_m)\mapsto (t,z_1,\cdots, z_m).$$
Then $e:=e_1\wedge\cdots\wedge e_m$ is a local frame of $\text{det}E$ over $U$, whose norm with respect to $\det h$ is given by
$$\|e(t)\|^2_{\det h}=\frac{1}{\mu_0(\phi_t(B_t))},$$
where $\mu_0$ is the Lebesgue measure on $\mc^m$.

By Corollary \ref{cor:volue vairation strict convex}, we know that $-\ln\mu_0(\phi_t(B_t))$ is a strictly plurisubharmonic function on $U$
provided that $h$ is strictly negatively curved.
Note that the curvature of $(\text{det}E,\text{det}h)$ on $U$ is given by $i\partial\bar\partial \ln\mu_0(\phi_t(B_t))$, so the curvature of the induced Hermitian metric $\text{det}h$ on $\text{det}E$ is strictly negative.
We thus get

\begin{cor}[=Corollary \ref{cor: detE strict negative}]
	Let $\pi:E\ra X$ be a holomorphic vector bundle over a complex manifold $X$ equipped with a smooth Finsler metric $h$.
	If $(E,h)$ is strictly negatively curved, then the curvature of the induced Hermitian metric $\text{det}h$ on $\text{det}E$ is strictly negative.
\end{cor}

\section{Deduce Theorem \ref{thm:ample bunlde positive metric} from Theorem \ref{thm:curvartue positive cicular domains} or Theorem \ref{thm:curvartue positive cicular domains from weight}} \label{sec:ample vector bundle}
In this section, we discuss the relation of Theorem \ref{thm:ample bunlde positive metric} with Theorem \ref{thm:curvartue positive cicular domains} or Theorem \ref{thm:curvartue positive cicular domains from weight}.
We show that  Theorem \ref{thm:ample bunlde positive metric} can be deduced
from Theorem \ref{thm:curvartue positive cicular domains} or Theorem \ref{thm:curvartue positive cicular domains from weight}.
For this consideration, the symmetric structure appearing in Theorem \ref{thm:curvartue positive cicular domains} or Theorem \ref{thm:curvartue positive cicular domains from weight} plays an indispensable role.

\subsection{Basic properties of ample vector bundles}\label{sec: ample defn}\label{subsec:ample v.b. and spsc}\
\par
This subsection recalls some well known basic knowledge about ample vector bundles.

Let $\pi:E\ra X$ be a holomorphic vector bundle over a compact complex manifold $X$.
For each $x\in X$, we denote by $E_x$ the fiber of $E$ over $x$ and denote by $E^*_x$ its dual.
Let $\mathbb P(E^*_x)$ be the projective space of $E^*_x$,
which is the space of one-dimensional complex linear subspaces of $E^*_x$ with the natural complex structure,
and let $\mathcal O_{\mathbb P(E^*_x)}(1)$ be the dual of the tautological line bundle over $\mathbb P(E^*_x)$.
Then
$$\mathbb P(E^*):=\cup_{x\in X}\mathbb P(E^*_x)$$
is a complex manifold that
can be naturally realized as a holomorphic fiber bundle over $X$ with $\mathbb P(E^*_x)$ as fibers,
and $$\mathcal{O}_{\mathbb P(E^*_x)}(1):=\cup_{x\in X}\mathcal O_{\mathbb P(E^*_x)}(1)$$
can be naturally realized as a holomorphic line bundle over $\mathbb P(E^*)$ whose restriction to $\mathbb P(E^*_x)$
is just $\mathcal O_{\mathbb P(E^*_x)}(1)$.
By definition, $E$ is called an \emph{ample vector bundle} if $\mathcal{O}_{\mathbb P(E*)}(1)$ is an ample line bundle over $\mathbb P(E^*)$.

We now assume that $E$ is ample.
Then there is a Hermitian metric $h$ on $\mathcal{O}_{\mathbb P(E^*)}(-1)$ whose curvature is negative.
Let $\rho:\mathcal{O}_{\mathbb P(E^*)}(-1)\ra \mr_{\geq 0}$ be the length function associated to $h$,
namely $\rho(v)=\sqrt{h(v,v)}$ for $v\in \mathcal{O}_{\mathbb P(E^*)}(-1)$.
Then $\rho$ is strictly plurisubharmonic on $\mathcal O_{\mathbb P(E^*)}(-1)\backslash Z_{\mathcal O_{\mathbb P(E^*)}(-1)}$,
where $Z_{\mathcal O_{\mathbb P(E^*)}(-1)}$ is the zero section of $\mathcal O_{\mathbb P(E^*)}(-1)$, viewed as a submanifold of
$\mathcal O_{\mathbb P(E^*)}(-1)$.

Note that $\mathcal{O}_{\mathbb P(E^*)}(-1)$ can be viewed as the blow up of $E^*$ along its zero section $Z_{E^*}$,
with $Z_{\mathcal O_{\mathbb P(E^*)}(-1)}$ as the exceptional divisor,
we can naturally identify $E^*\backslash Z_{E^*}$ with $\mathcal O_{\mathbb P(E^*)}(-1)\backslash Z_{\mathcal O_{\mathbb P(E^*)}(-1)}$.
Through this identification, we can view $\rho$ as a function on $E^*$, with $\rho|_{Z_{E^*}}\equiv 0$.

In conclusion, we get a plurisubharmonic function $\rho$ on $E^*$, which is
strictly plurisubharmonic on $E^*\backslash Z_{E^*}$ and invariant under the natural $S^1$ action on $E^*$.
In other words, $\rho$ is a smooth Finsler metric on $E^*$ whose curvature is strictly negative.


%

\subsection{Some linear algebra}\label{subsec:linear alg prep.}\
\par
We present some knowledge about linear algebra that is needed in the proof of  Theorem \ref{thm:ample bunlde positive metric}.

Let $V$ be a complex vector space with complex dimension $m$ and $V^*$ be its dual space.
Let $\p(V)$ be the set of all  polynomials on $V$, and $\p_k(V)$ be the space of homogeneous polynomials of degree $k$ on $V$. Then
$$\p(V) = \bigoplus_{k \geq 0} \p_k(V).$$
We have $\p_0(V)=\mc$ and $\p_1(V)=V^*$.
In an obvious manner, we can naturally identify $\p_k(V)$ with $S^kV$, the $k$-th symmetric product of $V^*$.

We realize the circle group $S^1$ as the space of complex numbers with unit norm.
Then $S^1$ acts on $V$ via scalar product.
This induces an action of $S^1$ on $\p(V)$ as follows:
$$\alpha \cdot f(v)= f(\alpha v),$$
where $ f \in \p(V)$, $v \in V$, $\alpha \in S^1$.
Then $\p_k(V)$ are character subspaces of $\p(V)$ associated to this action, namely, for $k\geq 0$ we have
$$\p_k(V)= \{f \in \p(V) \big| \alpha \cdot f =\alpha^k f, \forall \alpha \in S^1\}.$$

The cotangent bundle of $V$ is $$T^*V= V \times V^*.$$
It follows that the canonical bundle of $V$ is $$K_V =V \otimes \det{V^*},$$
where $\det{V^*}=\wedge^m V^*$.

 We now consider the coordinate representation of $K_V$.
 Let $u_1, \cdots , u_m$ be a basis of $V$ and $u_1^*, \cdots ,u_m^*$ be the associated dual basis of $V^*$.
 Then $u_1^* \wedge \cdots \wedge u_m^*$ is a basis of $\det{V^*}$.
Consider linear isomorphism:
$$V \longrightarrow \mc^n: z_1 u_1 + \cdots + z_m u_m \mapsto (z_1, \cdots , z_m),$$
then $u_1^* \wedge \cdots \wedge u_m^*$ corresponds to $dz_1 \wedge \cdots \wedge dz_m$, a basis of $\det{(\mc^m)^*}$.

Let $\Omega \subset V$ be an $S^1$ invariant domain containing 0,
then we have the following identification
\be\label{eq:canonical section to hol map}
 H^0(\Omega, K_\Omega)=\mathcal{O}(\Omega, \det{V^*}),
\ee
where $\mathcal{O}(\Omega, \det{V^*})$ is the space of holomorphic mappings from $\Omega$ to $\det V^*$.
The action of $S^1$ on $\mathcal{O}(\Omega, \det{V^*})$ is given as follows:
 $$\alpha \cdot f(x)=f(\alpha x), \ \alpha\in S^1.$$
Under coordinate form, if we identify $V$ with $\mc^m$ as above and view $\Omega$ as a domain in  $\mc^m$,
then we have the following identification
\be\label{eq:canonical section to function}
H^0 (\Omega,K_\Omega)\cong\{f(z_1, \cdots , z_m)dz_1 \wedge \cdots \wedge dz_m |f \in \mathcal{O}(\Omega)\},
\ee
and the action of $S^1$ on $H^0 (\Omega,K_\Omega)$ is realized as
$$\alpha \cdot (f(z_1, \cdots , z_m)dz_1 \wedge \cdots \wedge dz_m)=f(\alpha z_1, \cdots , \alpha z_m)dz_1 \wedge \cdots \wedge dz_m.$$
It is clear that the action of $S^1$ on the Hilbert space
$$A^2(\Omega)=\{f \in \operatorname{H}^0(\Omega,K_\Omega): ||f||  < +\infty\}$$
is unitary, where
$$||f||^2=\int_{\Omega} c_m f \wedge \bar{f},$$
with $c_m=\frac{i^{m^2}}{2^m}$ is set to make the form $c_m f \wedge \bar{f}$ real and nonnegative.
\par For any $k\geq 0,$ let
\be\label{eq:def of P'_k}
\p_k'(\Omega)=\{f \in A^2(\Omega) |\ \alpha \cdot f = \alpha^k f,\ \forall \alpha \in S^1\},
\ee
then
\be\label{eq:canonical section hom. poly.}
\p^{'}_k(\Omega)=\{f(z_1, \cdots , z_m)dz_1 \wedge \cdots \wedge dz_m|\ f \in \p_k(\mc^m)\}.
\ee
It follows that
\be\label{eq:canonical bdl and sym prod}
\p^{'}_k(\Omega)=\p_k(V) \otimes \det{V^*}=\s^kV^* \otimes \det{V^*}.
\ee

\subsection{The proof of Theorem \ref{thm:ample bunlde positive metric}}\
\par Let $\pi:E\ra X$ be an ample holomorphic vector bundle of rank $m$ over a compact complex manifold $X$ of dimension $n$.
Let $E^*$ be the dual bundle of $E$ and let $Z_{E^*}$ be the zero section of $E^*$, viewed naturally as a submanifold $E^*$.
From \S \ref{subsec:ample v.b. and spsc}, we know that $E^*$ admits a smooth Finsler metric $\rho:E^*\ra\mr_{\geq 0}$
whose curvature is strictly negative.

Let $\Omega=\{v\in E^*|\rho(v)\leq 1\}$, then $\Omega$ is an $S^1$ invariant bounded domain in $E^*$
whose boundary is strictly pseudoconvex.
As usual, we denote $\Omega\cap E^*_t$ by $\Omega_t$ for $t\in X$.
Note that $\Omega_t$ is an $S^1$-invariant domain in $E^*_t$ containing the origin.
By \eqref{eq:canonical section to hol map}, we can canonically identify $H^0(\Omega_t, K_{\Omega_t})$
with $\mathcal O(\Omega_t, \det E_t)$.
For $k\geq 0$, if we define $\p'_k(\Omega_t)$ as in \eqref{eq:def of P'_k},
we have $\p'_k(\Omega_t)=S^k E_t\otimes \det E_t$ from \eqref{eq:canonical bdl and sym prod}.

Let $\varphi$ be an $S^1$-invariant smooth plurisubharonic function defined on some neighborhood of the closure $\overline\Omega$ of $\Omega$ in $E^*$.
On $\p'_k(\Omega_t)$, we can define a Hermitian inner product $h_t$ by setting
$$\|f\|^2_{h_t}=\int_{\Omega_t}c_mf\wedge\bar f e^{-\varphi_t},\ f\in \p'_k(\Omega_t),$$
where $\varphi_t$ is the restriction of $\varphi$ on $\Omega_t$.
In this way, we get a Hermitian metric $h$ on $S^kE\otimes\det E$.
Our propose is to deduce from Theorem \ref{thm:curvartue positive cicular domains} or Theorem \ref{thm:curvartue positive cicular domains from weight}
that the curvature of the Hermitian vector bundle $(S^kE\otimes\det E, h)$ over $X$
is strictly positive in the sense of Nakano, for suitable choice of $\varphi$ (indeed for all such $\varphi$),
and hence get new proofs of Theorem \ref{thm:ample bunlde positive metric}.

The argument goes as follows.
Let $(U, t_1, \cdots , t_n)$ be a local coordinate on $X$,
and $e_1,\cdots , e_m$ be a holomorphic local frame of $E^*$ over $U$.
Then we get an isomorphism $\sigma:\pi^{-1}(U) \longrightarrow U \times \mc^m$ given by
$$(t, z_1 e_1 +\cdots +z_me_m) \mapsto (t_1,\cdots, t_n,z_1,\cdots,z_m),$$
where $\pi:E^*\ra X$ is the bundle map.
This isomorphism realizes $\Omega\cap \pi^{-1}(U)$ as a strictly pseudoconvex family of bounded domains over $U$
whose fibers $\sigma(\Omega_t)\subset\mc^m$ are circular domains containing the origin.
By \eqref{eq:canonical section hom. poly.}, for $t \in U$, via $\sigma$ we can identify $\p^{'}_k(\Omega_t)$ with the space
$$\{f(z_1,\cdots, z_m)dz_1\wedge\cdots dz_m|f\in \p_k(\mc^m)\},$$
with the Hermitian inner product $h_t$ given by
$$\|f(z_1,\cdots, z_m)dz_1\wedge\cdots dz_m\|^2_{h_t}=\int_{\sigma(\Omega_t)}|f|^2e^{-\varphi_t\circ\sigma^{-1}}d\lambda_z.$$
It follows from Theorem \ref{thm:curvartue positive cicular domains} that the curvature of $(S^kE\otimes\det E, h)$
is strictly positive in the sense of Nakano, and hence we get Theorem \ref{thm:ample bunlde positive metric}.

In a similar way, we can deduce Theorem \ref{thm:ample bunlde positive metric} from Theorem \ref{thm:curvartue positive cicular domains from weight}
by choosing
$\varphi=\max_{1/4, 1/4}\{1/4, \rho\}.$
(see Lemma \ref{lem: Regularization max functions} for the definition of the regularized maximum function).

	\end{document}